%% file: main.tex
\documentclass[10pt,a4paper]{article}
\usepackage[margin=2.5cm]{geometry}
\usepackage{amsmath,amssymb,amsthm,mathtools}
\usepackage{tikz-cd}
\usepackage[colorlinks=true,linkcolor=blue,citecolor=blue,urlcolor=blue]{hyperref}
\usepackage{cleveref}
\usepackage[titletoc]{appendix}

\usepackage[colorinlistoftodos,textsize=small]{todonotes}

\theoremstyle{plain}
\newtheorem{thm}{Theorem}[section]
\newtheorem{lem}[thm]{Lemma}
\newtheorem{prop}[thm]{Proposition}
\newtheorem{cor}[thm]{Corollary}

\newtheorem{mainthm}{Theorem}

\crefname{mainthm}{Theorem}{Theorems}

\theoremstyle{definition}
\newtheorem{defn}[thm]{Definition}
\newtheorem{rem}[thm]{Remark}
\newtheorem{expl}[thm]{Example}

\newtheorem{ass}{Hypothesis}

\makeatletter
\def\dt@alphabets#1{%
  \expandafter\newcommand\csname c#1\endcsname{\mathcal{#1}}%
  \expandafter\newcommand\csname b#1\endcsname{\mathbb{#1}}%
  \expandafter\newcommand\csname f#1\endcsname{\mathfrak{#1}}}
\@for\dt@letter:=A,B,C,D,E,F,G,H,I,J,K,L,M,N,O,P,Q,R,S,T,U,V,W,X,Y,Z\do{%
  \expandafter\dt@alphabets\dt@letter}
\makeatother

\DeclareMathOperator{\Aut}{Aut}

\DeclareMathOperator{\Fund}{Fund}

\DeclareMathOperator{\pr}{pr}

\DeclareMathOperator{\tor}{tor}

\title{Isogeny rigidity under algebraic group actions and finite translation orbits\\ \large (with an appendix by E.~Amerik)}
\author{Paolo Dolce \and Francesco Tropeano}
\date{}

\newcommand{\Addresses}{{
  \bigskip
  \bigskip
  \footnotesize
}
  E.~Amerik, \textsc{Universit\'e Paris-Saclay, Laboratoire de Mathématiques d’Orsay, France} and   \textsc{National Research University Higher School of Economics, Laboratory of Algebraic Geometry and its Applications, Russia }\par\nopagebreak
  \textit{E-mail address}: \texttt{ekaterina.amerik@gmail.com
}
  \bigskip

  P.~Dolce, \textsc{Westlake University, Institute for theoretical sciences, China}\par\nopagebreak
  \textit{E-mail address}: \texttt{dolce@westlake.edu.cn}

  \bigskip
  
F.~Tropeano, \textsc{Westlake University, Institute for theoretical sciences, China}\par\nopagebreak
  \textit{E-mail address}: \texttt{francesco.tropeano@westlake.edu.cn}
}

\begin{document}

\maketitle
\begin{abstract}
We study finite orbits generated by fibrewise translations on varieties carrying two abelian fibrations. Under natural dominance and generation hypotheses, non-isotriviality of the first family forces points with finite orbit to lie over proper closed subsets of both bases.  This extends results of Corvaja, Tsimerman and Zannier to complex families of arbitrary, possibly different relative dimensions, including torsors. The proof combines relative Manin--Mumford with a general isogeny rigidity theorem for algebraic group actions.
\end{abstract}

\makeatletter
\@starttoc{toc}
\makeatother

\section{Introduction}
\input{Introduction}

\section{Notation and conventions}
\input{Notations}

\section{Local monodromy around divisors and  rigidity under group actions}
\input{Rigidity.tex}

\section{Translation orbits on double families of abelian varieties}\label{sec:graph}
\input{DoubleFamilies.tex}

\clearpage
\begin{appendices}

\section{Construction of double abelian fibrations in the IHS case}
\label{sec:app}

\begin{center}
    \large E.~Amerik
\end{center}
\medskip

\input{example-paolo-francesco.tex}

\end{appendices}


\bibliographystyle{amsplain}
\bibliography{biblio}
\Addresses
\end{document}

%% file: Introduction.tex
\theoremstyle{plain}
\newtheorem*{introRigidity}{Theorem}
\newtheorem*{introTranslation}{Corollary}
\newtheorem*{introA}{Theorem A}
\newtheorem*{introB}{Theorem B}
\newtheorem*{introC}{Theorem C}
\newtheorem*{introD}{Theorem D}
\theoremstyle{definition}

The distribution of finite orbits is a basic problem in the dynamics of groups of automorphisms of algebraic varieties. A central question is whether the points with finite orbit are contained in a proper closed subset, and what geometry allows them to be Zariski dense. For projective surfaces, Cantat and Dujardin~\cite{CD} relate this density to Kummer geometry: when the surface and the group are defined over a number field and the group contains both parabolic and loxodromic automorphisms, density of finite orbits forces the pair to be a Kummer group.

A geometrically interesting case is provided by an abelian scheme $\mathcal B/S$ of relative dimension $g$ and a section $\sigma$. A point $p\in\mathcal B_s$ has finite orbit under the fibrewise translation $t_\sigma$ precisely when $\sigma(s)$ is torsion. The dynamical question is therefore exactly the problem of determining whether the torsion specializations of $\sigma$ are Zariski dense. This is a known problem in unlikely intersections called Relative Manin--Mumford conjecture, which has been recently settled  by Gao and Habegger in~\cite{GH}. Their stronger density criterion shows that, if the multiples of $\sigma$ are Zariski dense in $\mathcal B$, its torsion specializations are Zariski dense in $S$ if and only if its Betti map has rank $2g$. In particular, they lie in a proper closed subset when $\dim S<g$. The boundary case $\dim S=g$ is especially relevant here: if the Betti map has maximal rank, torsion specializations are analytically dense, so no confinement is possible for the single translation.

We study varieties carrying two abelian fibrations, with translations $t_1$ and $t_2$ induced by sections $\sigma_1$ and $\sigma_2$. In the elliptic case, Corvaja, Tsimerman and Zannier~\cite{CTZ} study this problem for surfaces with two non-isotrivial elliptic fibrations $f_i\colon X\to\mathbb P^1$ and non-torsion sections $\sigma_1$ and $\sigma_2$, with all data defined over a number field. Assuming that $(f_1,f_2)$ is generically finite, they prove that the points $p$ for which
\[
\bigl\{t_1^n\circ t_2^r(p): n,r\geq0,\ \text{the composition is defined at }p\bigr\}
\]
is finite lie on finitely many fibres of $f_2$. Thus their confinement theorem requires finiteness only for these ordered words, a weaker condition than finiteness of the full $\langle t_1,t_2\rangle$-orbit. Their approach connects this condition with simultaneous torsion specializations, combining Pila--Wilkie counting, height estimates, and functional transcendence results of Ax--Schanuel type.

We generalize this confinement result in several directions. All the data are defined over $\mathbb C$; the two fibrations may have arbitrary, possibly different relative dimensions over higher dimensional bases; their smooth parts may be torsors under abelian schemes, so global zero sections are not required. We assume that the joint map to the two bases is dominant and that each translating section generates the geometric generic fibre of its associated abelian scheme. Under these hypotheses, if the first associated abelian scheme is non-isotrivial, points with finite ordered orbit are confined over a proper closed subset of the second base.

The main difficulty occurs precisely at the boundary left open by the individual torsion conditions. Relative Manin--Mumford first shows that density of the finite-orbit locus forces the two base dimensions and the two relative dimensions to coincide, with maximal Betti rank for both sections. In this case, points that are periodic for each translation separately can be analytically dense in the common total space. We overcome this difficulty by combining relative Manin--Mumford with an isogeny rigidity theorem: the former turns density of finite orbits into algebraic relations between the families, while the latter forces the first family to be isotrivial. This rigidity theorem is a main result in its own right, formulated for arbitrary algebraic group actions on smooth bases, independently of the double-family setting.

\subsection{Main results on translation orbits}

We now specify the double-family setting. Let $\mathcal A'$ be a smooth projective variety, and let
\[
f_i\colon\mathcal A'\dashrightarrow S_i',\qquad i=1,2,
\]
be dominant rational maps to projective varieties.  Assume that each $f_i$ has a smooth part with fibres isomorphic as varieties to abelian varieties of dimension $g_i\geq1$. Write $d_i:=\dim S_i$, so that $\dim\mathcal A'=g_1+d_1=g_2+d_2$. The family $\mathcal A_i \to S_i$ is canonically a torsor under its abelian scheme of translations $\mathcal B_i/S_i$. A section $\sigma_i\in\mathcal B_i(S_i)$ determines an automorphism
\[
t_i\colon\mathcal A_i\longrightarrow\mathcal A_i, \qquad p\longmapsto p+\sigma_i(f_i(p)),
\]
and hence a birational transformation of $\mathcal A'$. The addition denotes the torsor action and requires no choice of origins. The two hypotheses used for global confinement are:
\begin{enumerate}
    \item[\textup{(H1)}] The rational map $\Psi:=(f_1,f_2)\colon\mathcal A'\dashrightarrow S_1'\times S_2'$ is dominant, i.e. the fibrations are distinct.

    \item[\textup{(H2)}] For $i=1,2$,  the cyclic group $\mathbb Z\sigma_i$ is Zariski dense in $\mathcal B_i$.
\end{enumerate}
We have the group $\Gamma:=\langle t_1,t_2\rangle\subseteq\operatorname{Bir}(\mathcal A')$ and we study the subsets 
\[
O_{12}:=\{t_1^n\circ t_2^r:n,r\in\mathbb N\}, \qquad \mathcal T :=\{t_1^n:n\in\mathbb N\} \cup\{t_2^r:r\in\mathbb N\} \cup\{t_1^n\circ t_2:n\in\mathbb N\}.
\]
 For a collection $\Theta$ of words in $\Gamma$ define
\[
\Theta(p):=\{w(p):w\in\Theta\text{ is evaluable at }p\}, \qquad \mathfrak F(\Theta):=\{p\in\mathcal A':\#\Theta(p)<\infty\}.
\]
Thus $\mathfrak F(\Gamma)\subseteq\mathfrak F(O_{12}) \subseteq\mathfrak F(\mathcal T)$. Notice that $\mathfrak F(\mathcal T)$ is contained in the intersection of the fibers over torsion values of $\sigma_1$ and $\sigma_2$, which might be Zariski dense in $\mathcal A'$. However, using some auxiliary families, the points of $\mathfrak F(\mathcal T)$ can be precisely characterized in the following way: on the open subset $U_{12}:=\{p\in\mathcal A_1\cap\mathcal A_2:t_2(p)\in\mathcal A_1\}$, we define the abelian schemes
\[
\mathcal P_{12} := f_1^*\mathcal B_1\times_{U_{12}}f_2^*\mathcal B_2, \qquad \mathcal P_3 :=(f_1\circ t_2)^*\mathcal B_1,
\]
with sections $\xi_{12}:=(f_1^*\sigma_1,f_2^*\sigma_2)$ and $\xi_3 :=(f_1\circ t_2)^*\sigma_1$, respectively. Then, for $p\in U_{12}$, one has
\[
p\in\mathfrak F(\mathcal T) \quad\Longleftrightarrow\quad \sigma_1(f_1(p)),\ \sigma_2(f_2(p)),\ \sigma_1(f_1(t_2(p)))\text{ are all torsion}.
\]
We say that the \emph{graph relation} holds if there exist a dense open subset $U'\subseteq U_{12}$, an integer $N\geq1$, and a homomorphism of abelian schemes $\theta\colon\mathcal P_{12}|_{U'}\longrightarrow\mathcal P_3|_{U'}$ such that
\[
[N]\xi_3|_{U'}=\theta(\xi_{12}|_{U'}).
\]
Our first result shows that satisfying the graph relation is equivalent to the density of $\mathfrak F(\mathcal T)$.

\begin{introA}
In the double-family setting above, assume \textup{(H1)} and \textup{(H2)}. The following conditions are equivalent:
\begin{enumerate}
\item[$(i)$] $\mathfrak F(O_{12})$ is Zariski dense in $\mathcal A'$.

\item[$(ii)$] $\mathfrak F(\mathcal T)$ is Zariski dense in $\mathcal A'$.

\item[$(iii)$] There exists an integer $g\geq1$ such that
\begin{equation}\label{intro:balanced}
    g_1=g_2=d_1=d_2=g, \qquad \operatorname{rk}\beta_{\sigma_1} =\operatorname{rk}\beta_{\sigma_2}=2g,
\end{equation}
and the graph relation holds.
\end{enumerate}
\end{introA}
The following example shows an explicit construction where the above three equivalent conditions hold. 
\begin{expl}\label{expl:constantExa}
Let $B$ be a positive-dimensional complex abelian variety, take $\mathcal A'=B\times B$ with its two projections, and consider the translations
\[
t_1(a,b)=(a,b+a),\qquad t_2(a,b)=(a+b,b),
\]
induced by the diagonal sections $\sigma_1(a)=(a,a)$ and $\sigma_2(b)=(b,b)$, respectively. One can verify that \textup{(H1)}, \textup{(H2)} and \eqref{intro:balanced} hold. The graph relation holds over all of $\mathcal A'$ with $N=1$. Indeed, the natural trivializations identify
\[
\mathcal P_{12}\simeq(B\times B)\times\mathcal A',
\qquad
\mathcal P_3\simeq B\times\mathcal A'.
\]
In fibre coordinates, the sections are $\xi_{12}(a,b)=(a,b)$ and $\xi_3(a,b)=a+b$. Thus the homomorphism of abelian schemes
\[
\theta\colon(B\times B)\times\mathcal A'
\longrightarrow B\times\mathcal A',
\qquad
((u,v),p)\longmapsto(u+v,p),
\]
satisfies $\xi_3=\theta(\xi_{12})$. At the same time, one can check directly that 
\[
\bigcup_{N\geq1}B[N]\times B[N]
\subseteq\mathfrak F(\Gamma)
\]
is Zariski dense in $B\times B$.
\end{expl}
Our second main result says that the graph relation cannot hold if the first family is non-isotrivial: 

\begin{introB}
For a double family with translation vectors $\sigma_1,\sigma_2$ as above, the following assertions hold:
\begin{enumerate}
\item[$(i)$] If $\mathcal B_1/S_1$ is constant, then the graph relation holds.

\item[$(ii)$] Assume \textup{(H1)} and \textup{(H2)}. If the graph relation holds, then $\mathcal B_1/S_1$ is isotrivial.
\end{enumerate}
\end{introB}

Example~\ref{expl:constantExa} provides a constant double family for which the graph relation holds and 
$\mathfrak F(\Gamma)$ is Zariski dense.
However, the constancy assumption in Theorem~B$(i)$ cannot be
weakened to isotriviality.
In Section~\ref{subsec:isoDensity}, we construct two double elliptic
families satisfying \textup{(H1)}, \textup{(H2)} and
\eqref{intro:balanced}, with both associated elliptic schemes
isotrivial but nonconstant. In the first example, the graph
relation holds and $\mathfrak F(\Gamma)$ is Zariski dense; in the second, the graph relation fails and $\mathfrak F(\mathcal T)$ is not Zariski dense.

To state confinement also at indeterminacy points, let $\Gamma_{f_i}\subseteq\mathcal A'\times S_i'$ be the closure of the graph of $f_i$, and define, for a closed subset $Z\subseteq S_i'$,
\[
f_i^{-1}[Z]
:=\operatorname{pr}_{\mathcal A'}
\bigl(\Gamma_{f_i}\cap(\mathcal A'\times Z)\bigr).
\]
Since $S_i'$ is projective, this is a closed subset of $\mathcal A'$; on the domain of definition of $f_i$, it agrees with the ordinary inverse image of $Z$. Combining Theorem~A and~B we then get:

\begin{introC}
In the double-family setting above, assume \textup{(H1)} and \textup{(H2)}. Suppose that condition \eqref{intro:balanced} fails or that the graph relation does not hold. Then:
\begin{enumerate}
\item[$(i)$] $\mathfrak F(\mathcal T)$ is not Zariski dense in $\mathcal A'$.

\item[$(ii)$] There exists a proper closed subset $Z_2\subsetneq S_2$ such that $\mathfrak F(O_{12})\cap\mathcal A_2\subseteq f_2^{-1}(Z_2)$.

\item[$(iii)$] There exist proper closed subsets $W_i\subsetneq S_i$, for $i=1,2$, such that $\mathfrak F(\Gamma)\cap\mathcal A_i\subseteq f_i^{-1}(W_i)$.
\end{enumerate}
Moreover, there exist proper closed subsets $Z_2'\subsetneq S_2'$ and $W_i'\subsetneq S_i'$ for $i=1,2$ satisfying $Z_2\subseteq Z_2'\cap S_2$ and $W_i\subseteq W_i'\cap S_i$ such that
\[
\mathfrak F(O_{12})\subseteq f_2^{-1}[Z_2'],
\qquad
\mathfrak F(\Gamma)
\subseteq f_1^{-1}[W_1']\cap f_2^{-1}[W_2'].
\]
In particular, all these conclusions hold if $\mathcal B_1/S_1$ is non-isotrivial.
\end{introC}

For $g_1=g_2=d_1=d_2=1$, the sets $Z_2'$ and $W_i'$ are finite. Thus Theorem C confines finite ordered orbits to finitely many fibres of $f_2$. If $\Psi$ is a finite morphism, then
\[
f_1^{-1}(W_1')\cap f_2^{-1}(W_2')
=\Psi^{-1}(W_1'\times W_2')
\]
is finite, giving finiteness of the full group-orbit locus. Thus we recover the confinement of \cite[Theorem~1.1]{CTZ} and its consequence for finite group orbits in \cite[Corollary~1.2]{CTZ}. Moreover, part~$(i)$ recovers the confinement to finitely many curves in \cite[Theorem~1.4]{CTZ}: after the birational change of variable $p=t_1(x)$, their three torsion conditions are precisely those defining $\fF(\cT)$ on the common smooth domain. Here, all these conclusions hold without a number-field assumption on the data and assuming non-isotriviality only for $\cB_1/S_1$.

The final theorem on the translation orbits concerns an individual fibre on which finite ordered orbits are dense. It uses weaker hypotheses than the global confinement theorem. For a torsor $A$ under an abelian variety $B$ and an abelian variety $G/\mathbb C(A)$, write
\[
\operatorname{Isog}(G)
:=\{x\in B(\mathbb C):T_x^*G\text{ is }\mathbb C(A)\text{-isogenous to }G\},
\]
and let $\operatorname{Bad}(G)$ be the finite set of prime divisors of $A$ at which $G$ has bad reduction.

\begin{introD}
Assume that $g_1=g_2=g\geq1$ and that $\mathbb Z\sigma_1$ is Zariski dense in $\cB_1$. Fix $b \in S_2(\mathbb C)$ such that the restriction $\phi_b:=f_1|_{\cA_{2,b}}\colon \cA_{2,b}\dashrightarrow S_1$ is dominant. Put $K_b:=\mathbb C(\cA_{2,b})$ and $G_b:=(\cB_{1,\eta_1})_{K_b}$, where $\eta_1$ is the generic point of $S_1$, and the extension $\mathbb C(S_1)\hookrightarrow K_b$ is induced by $\phi_b$. If $\fF(O_{12})\cap \cA_{2,b}$ is Zariski dense in $\cA_{2,b}$, then
\[
\sigma_2(b)\in \mathcal B_{2,b}(\mathbb C)_{\tor} \cap\operatorname{Isog}(G_b).
\]
In particular, if $\mathcal B_{2,b}$ is simple and $G_b$ is non-isotrivial, then
\[
\operatorname{ord}(\sigma_2(b)) \leq |\operatorname{Isog}(G_b)| \leq |\operatorname{Bad}(G_b)| \min_{D\in\operatorname{Bad}(G_b)} \frac{D^g}{g!},
\]
where $D^g$ is the top self-intersection number of the reduced prime divisor $D$ on $\mathcal A_{2,b}$.
\end{introD}
Note that this result imposes more than torsion of $\sigma_2(b)$: the corresponding translation must preserve the isogeny class of the first family restricted to $\mathcal A_{2,b}$. When this translation--isogeny group is finite, its order controls the permitted torsion orders. The estimate depends on the bad divisors of the individual restricted family.

In an appendix to this paper, E. Amerik constructs explicit examples satisfying our hypotheses in every even dimension. Starting from a quartic K3 surface $S$ containing two skew lines, she obtains two Lagrangian fibrations $S^{[g]}\to\mathbb P^g$ on its Hilbert scheme of length-$g$ subschemes. Over their smooth loci, both fibrations are abelian schemes without a fixed part, each with a zero section and an additional non-torsion section whose multiples are Zariski dense. These constructions provide non-isotrivial examples to which our confinement theorems apply in every relative dimension $g\geq1$.

\subsection{Isogeny rigidity under group actions}

The geometric ingredient in the confinement theorems is a rigidity principle for families whose isogeny class is preserved by automorphisms of the base. Together with relative Manin--Mumford, it rules out the graph relation when the first family is non-isotrivial. It also gives a general description of isogeny symmetries, independently of the double-family setting. We work with a pair $(X,H)$, where $X$ is a smooth irreducible quasi-projective complex variety and $H$ is a complex algebraic group acting algebraically on $X$.

Put $K=\mathbb C(X)$, let $G/K$ be an abelian variety, and choose an abelian scheme model $\pi\colon\mathcal G\to U$ over a dense open subset $U\subseteq X$. Our starting point is the finite set
\[
\operatorname{Bad}(G)
:=\{D\subseteq X:D\text{ is a prime divisor at which }G
\text{ has bad reduction}\},
\]
which generalizes the places of bad reduction over a curve to divisors on a base of arbitrary dimension. It is independent of the chosen model and invariant under $K$-isogeny. Thus an automorphism preserving the isogeny class must also preserve this set of divisors.

The elliptic case illustrates how this restriction leads to rigidity. Suppose that $E$ is an elliptic curve and that the isogeny class of an elliptic curve $G/\mathbb C(E)$ is preserved by a non-torsion translation of $E$. Its finite set of bad places must be invariant under this translation, hence empty. The resulting elliptic scheme over $E$ has constant $j$-invariant and is therefore isotrivial. In higher dimension, bad divisors replace bad points, and period maps replace the $j$-invariant. Since a divisor may now be preserved by a positive-dimensional group, the corresponding conclusion is isotriviality along its orbits.

To pass from bad reduction to period maps, we use local monodromy. If $M_D$ denotes the monodromy around a general point of a prime divisor $D$, then
\[
D\notin\operatorname{Bad}(G)\quad\Longleftrightarrow\quad M_D=1.
\]
This triviality allows the integral local system $\mathbb V:=R^1\pi_*\mathbb Z$ to extend uniquely to
\[
W:=X\setminus\bigcup_{D\in\operatorname{Bad}(G)}D
\]
by Lemma~\ref{lem:extension-local-system}. The next step is to turn extension into isotriviality. Proposition~\ref{prop:Liouville-general} shows that if $\mathbb V$ extends over $X$ and every bounded holomorphic function on its analytic universal cover is constant, then $G$ is isotrivial. The proof applies this Liouville property to a bounded realization of the period map.

We now bring the group action into this argument. Define
\[
\begin{aligned}
\operatorname{Isom}_H(G)
&:=\{h\in H(\mathbb C):h^*G\simeq_K G\},\\
\operatorname{Isog}_H(G)
&:=\{h\in H(\mathbb C):h^*G\sim_K G\},
\end{aligned}
\qquad
\Sigma_G:=\operatorname{Stab}_H(\operatorname{Bad}(G)).
\]
Here $\sim_K$ denotes $K$-isogeny, $h^*G$ is the pullback under the automorphism induced by $h$, and $\Sigma_G$ is the setwise stabilizer. The universal cover of every connected complex algebraic group has the Liouville property. Applying the preceding criterion along orbit maps from $\Sigma_G^0$ therefore gives isotriviality along its orbits meeting $U$. Theorem~\ref{rig:thm:invariance} strengthens this conclusion: both $\operatorname{Isom}_H(G)$ and $\operatorname{Isog}_H(G)$ are closed algebraic subgroups of $H$, and, identifying them with their complex points,
\[
\Sigma_G^0\subseteq\operatorname{Isom}_H(G)
\subseteq\operatorname{Isog}_H(G)\subseteq\Sigma_G,
\qquad
\operatorname{Isom}_H(G)^0
=\operatorname{Isog}_H(G)^0
=\Sigma_G^0.
\]
Thus the connected isogeny symmetries are determined by the bad divisors. If their orbits are dense, isotriviality along orbits becomes isotriviality of the whole family. This gives the key rigidity statement, Corollary~\ref{thm:translation}.

\begin{introTranslation}
If $\operatorname{Isog}_H(G)$ has a Zariski-dense orbit in $X$, then $G$ is isotrivial. In particular, if $h\in H(\mathbb C)$ is such that $\overline{\langle h\rangle}^{\mathrm{Zar}}$ has a Zariski-dense orbit in $X$ and $h^k\in\operatorname{Isog}_H(G)$ for some $k\in\mathbb Z\setminus\{0\}$, then $G$ is isotrivial.
\end{introTranslation}

This corollary supplies the step from the graph relation to confinement. In Theorem A, relative Manin--Mumford turns density of the simultaneous torsion values into the graph relation. For Theorem B, we restrict this relation to a very general fibre of the second fibration, where the second translation generates a Zariski-dense subgroup. We then quotient the pulled-back first family by its isotrivial part. The graph relation yields an isogeny under translation on the abelian subvariety generated by the induced section. The corollary makes this subvariety isotrivial, whereas the quotient has no nonzero isotrivial abelian subvariety. The induced section is therefore torsion, and the generation hypothesis forces the quotient to vanish. Hence the first family is isotrivial. For a non-isotrivial first family, this contradiction excludes the graph relation and leads to the confinement statements of Theorem C.

Theorem D uses the same ideas to study a fibre of the second fibration on which finite ordered orbits are Zariski dense. Under its hypotheses, relative Manin--Mumford shows that the translation vector on this fibre is torsion and that its translation preserves the isogeny class of the restricted first family. The vector therefore belongs to the corresponding translation--isogeny group. If this group is finite, the order of the vector divides its cardinality. Bounding that cardinality in terms of the bad divisors of the restricted family gives the torsion-order estimate in Theorem D.

\subsection{Comparison with results in Hyperk\"ahler dynamics}
Amerik and Cantat~\cite[\S6.2]{AC} study a similar finite-orbit problem for an irreducible hyperk\"ahler manifold $M$ of dimension $2g$ carrying two parabolic automorphisms with distinct Lagrangian fibrations. After taking suitable positive powers, these automorphisms act by fibrewise translations $t_1,t_2$. For $\Gamma_{k,\ell}:=\langle t_1^k,t_2^\ell\rangle$,
they prove that, for every $\varepsilon>0$, the finite-orbit locus $\mathfrak F(\Gamma_{k,\ell})$ is $\varepsilon$-dense in $M$ for all sufficiently large $k,\ell$. Thus every point of $M$ lies within distance $\varepsilon$ of a finite orbit once both powers are sufficiently large, for a fixed metric inducing the Euclidean topology.

Their argument uses maximal Betti rank: the loci where the translation vectors have orders dividing $k$ and $\ell$ become increasingly dense in the respective bases. Since the joint fibration is dominant, these conditions can be satisfied simultaneously. In fact, the points they construct are common fixed points of $t_1^k$ and $t_2^\ell$, so already the fixed-point loci form arbitrarily fine nets in $M$ as the powers vary.

Our results describe the algebraic constraints on these loci. The bases and fibres have dimension $g$, and maximal Betti rank implies the generation condition \textup{(H2)}; together with dominance of the joint fibration, this places their setting within our hypotheses. Suppose that one associated abelian scheme is non-isotrivial, and label it as the first. For each fixed $k,\ell\geq1$, we apply Theorem C to the translations $t_1^k,t_2^\ell$, whose sections $k\sigma_1,\ell\sigma_2$ still satisfy \textup{(H2)}. The theorem confines points with finite ordered orbit under these translations over a proper closed subset of the second base, and points with finite $\Gamma_{k,\ell}$-orbit over proper closed subsets of both bases. In particular,
\[
\mathfrak F(\Gamma_{k,\ell})\subseteq Z_{k,\ell}\subsetneq M
\]
for a proper Zariski-closed subset $Z_{k,\ell}$. The two conclusions distinguish a fixed dynamical system from the family obtained by varying the generating powers. For each fixed pair, the finite-orbit locus is contained in a proper algebraic subset and is therefore not Euclidean dense. Nevertheless, Amerik--Cantat show that for any prescribed positive scale, every sufficiently large pair has finite orbits within that distance of every point of $M$. The exceptional loci are thus algebraically confined for each group but collectively Euclidean dense. This also shows that our confinement cannot be uniform in the powers: the sets $Z_{k,\ell}$ cannot all be contained in one proper closed subset of $M$.

\paragraph{Organization of the paper.}
Subsections~\ref{subsec:BadReduction} and~\ref{subsec:translation-groups-K} develop divisorial monodromy and isogeny rigidity, providing the geometric tools for the main proofs. Section~\ref{sec:graph} develops the double-family setting. At the end of Subsection~\ref{subsec:4.2}, we deduce Theorem~\ref{thm:C} from Theorems~\ref{thm:A} and~\ref{thm:B}. The proofs of these two theorems are then given, in that order, in Subsection~\ref{subsec:4.3}. Finally, Subsection~\ref{subsec:exceptional-fibres} is devoted to the proof of Theorem~\ref{thm:D}.

\paragraph{Acknowledgements}
The authors are grateful to the anonymous referee for valuable comments and for identifying serious issues in an earlier version of this paper. In particular, the referee brought Example~\ref{expl:constantExa} to their attention, revealing a fundamental obstruction to the original proof strategy. These observations prompted a substantial revision of both the approach and the results of the paper.

%% file: Notations.tex
We fix the following notation and conventions: 
\begin{itemize}

\item[$(1)$]
Unless another ground field is specified, all \emph{varieties} are integral separated schemes of finite type over $\bC$, and all morphisms of varieties are understood to be over $\bC$. Algebraic groups are understood to be separated group schemes of finite type over $\bC$ and need not be connected. Unless otherwise stated, a closed subset is Zariski closed, and pointwise statements about a $\bC$-variety $X$ concern its closed points, identified with $X(\bC)$. Generic points and fibres over function fields will be indicated explicitly. For a field extension $L/K$ and a variety $X$ over $K$, we write
\[
X_L:=X\times_{\operatorname{Spec}K}\operatorname{Spec}L.
\]
We use the same notation for a rational point and its image after extending the ground field. For a divisor $D$ on a smooth projective variety of dimension $m\geq1$, we write $D^m$ for its \emph{top self-intersection number}.

\item[$(2)$]
For a rational map $f\colon X\dashrightarrow Y$, we write $\Fund(f)$ for its \emph{fundamental locus}, namely the closed subset where $f$ is not defined, and $\Gamma_f\subseteq X\times Y$ for the closure of its graph. For a subset $Z\subseteq Y$, we define its \emph{set-theoretic total transform} by
\[
f^{-1}[Z] := \pr_X\bigl(\Gamma_f\cap\pr_Y^{-1}(Z)\bigr),
\]
where $\pr_X$ and $\pr_Y$ are the projections from $X\times Y$. On $X\setminus\Fund(f)$, this agrees with the ordinary inverse image of $Z$. Whenever an ordinary inverse image under a rational map is used, it is taken on the domain of definition of the map. If $Y$ is proper and $Z\subseteq Y$ is closed, then $f^{-1}[Z]$ is closed in $X$.

\item[$(3)$]
We write $\cB/S$ for an abelian scheme
$\pi\colon\cB\to S$ when the structure morphism is understood.
Isomorphisms of abelian varieties or abelian schemes are
understood to preserve the group structures.
We denote the zero section by $0_\cB$ and multiplication
by $n\in\bZ$ by $[n]\colon\cB\to\cB$. For $N\geq1$, we write
\[
\cB[N]:=\ker([N]\colon\cB\to\cB),
\qquad
\cB_{\tor}:=\bigcup_{N\geq1}\cB[N](\bC)
\]
for the \emph{$N$-torsion subgroup scheme} and the set of
\emph{fibrewise torsion points}, respectively.

We denote by $\cB(S):=\operatorname{Hom}_S(S,\cB)$ the
abelian group of sections, with fibrewise addition.
For $\sigma\in\cB(S)$, we write
\[
n\sigma:=[n]\circ\sigma,
\qquad
\bZ\sigma:=\bigcup_{n\in\bZ}(n\sigma)(S).
\]
The section $\sigma$ is \emph{torsion} if
$N\sigma=0_\cB$ for some $N\geq1$, equivalently if
$\sigma(S)\subseteq\cB[N]$ for some $N\geq1$; otherwise
it is \emph{non-torsion}. A \emph{torsion value} of
$\sigma$ is a point $s\in S(\bC)$ such that $\sigma(s)$
is torsion in $\cB_s$, and we denote the set of torsion
values by $T_\sigma:=\sigma^{-1}(\cB_{\tor})$.

Let $\eta$ be the generic point of $S$, put
$K:=\bC(S)$, and write $B:=\cB_\eta$ for the generic fibre.
For an algebraic closure $\overline K$ of $K$, let
$\overline\eta\colon\operatorname{Spec}\overline K\to S$
be the corresponding geometric generic point.
We identify $\cB_{\overline\eta}=B_{\overline K}$.
For a section $\sigma\in\cB(S)$, we denote its generic
and geometric generic values by
$\sigma(\eta)\in B(K)$ and
$\sigma(\overline\eta)\in B(\overline K)$, respectively.

\item[$(4)$] An abelian scheme $\cB/S$ over a $\bC$-variety $S$ is \emph{constant} if there exists an abelian variety $B_0$ over $\bC$ and an isomorphism $\cB\simeq B_0\times_{\bC}S$ of abelian schemes over $S$. An abelian variety $A$ over a field extension $K/\bC$ is \emph{isotrivial} if there exists an abelian variety $A_0$ over $\bC$ such that $A_{\overline K}\simeq(A_0)_{\overline K}$ as abelian varieties. An abelian scheme over an integral $\bC$-variety is \emph{isotrivial} if its generic fibre is isotrivial. More explicitly, let $\cB/S$ be an abelian scheme of relative dimension $g$ over a smooth $\bC$-variety $S$, and fix a polarization $\lambda$ of degree $\delta^2$. Write $K:=\bC(S)$ and $B:=\cB_\eta$. Let $\mu_\lambda\colon S\longrightarrow \mathbb A_{g,\delta}$ be its classifying morphism to the coarse moduli space of polarized abelian varieties of dimension $g$ and polarization degree $\delta^2$. Then the following conditions are equivalent and characterize isotriviality:
\begin{itemize}
    \item[$(i)$]
    The fibres $(\cB_s,\lambda_s)$, for $s\in S(\bC)$, are all isomorphic as polarized abelian varieties.

    \item[$(ii)$]
    The morphism $\mu_\lambda$ is constant, that is, it factors through a point of $\mathbb A_{g,\delta}(\bC)$.

    \item[$(iii)$]
    There exists an abelian variety $B_0$ over $\bC$ such that $B_{\overline K}\simeq(B_0)_{\overline K}$ as abelian varieties.

    \item[$(iv)$]
    There exist a finite \'etale surjective morphism $S'\to S$ and a polarized abelian variety $(B_0,\lambda_0)$ over $\bC$ such that
    \[
    (\cB,\lambda)\times_S S'
    \simeq
    (B_0,\lambda_0)\times_{\bC}S'
    \]
    as polarized abelian schemes over $S'$.
\end{itemize}
The equivalence with $(iii)$ shows that this notion
does not depend on the chosen polarization.
For the coarse moduli space, see
\cite[Theorem~20.13]{MilneAV1986}.
The finite \'etale criterion follows from the polarized
isomorphism scheme; see \cite[Theorem~2.1]{Sekiguchi1982}
and \cite[Tag~0GS9]{stacks-project}.

\item[$(5)$]
For an abelian variety $A$ over a field extension $K/\bC$,
we write $\operatorname{Tr}_{K/\bC}(A)$ for its
$K/\bC$-trace and
\[
\tau_A\colon
\bigl(\operatorname{Tr}_{K/\bC}(A)\bigr)_K
\longrightarrow A
\]
for the canonical trace homomorphism.
This homomorphism is universal for homomorphisms
from constant abelian varieties over $K$ to $A$.
We call its image the \emph{fixed part} of $A$ over $K$.

We denote by $I(A)\subseteq A$ the largest isotrivial
abelian subvariety of $A$ defined over $K$, and call
it the \emph{isotrivial part} of $A$.
For a suitable finite extension $L/K$, one has
\[
I(A)_L=\operatorname{im}(\tau_{A_L}).
\]
Thus the isotrivial part becomes the fixed part after
a suitable finite extension of the ground field.

\end{itemize}

%% file: Rigidity.tex
For an abelian variety over the function field of a complex curve, the places of bad reduction and their local monodromy provide information about the global geometry of the family. We study this relation for an abelian variety $G$ over $K=\mathbb C(V)$, where $V$ is a smooth complex variety of arbitrary dimension. In this setting, we consider the boundary divisors of an abelian scheme with generic fibre $G$ and describe their local monodromy. We show that triviality of this monodromy is precisely the condition for extending the associated integral local system across the boundary. The finite set $\operatorname{Bad}(G)$ of divisors with nontrivial local monodromy is independent of the chosen abelian scheme model. Under a suitable hypothesis on the universal cover of $V$, its emptiness forces $G$ to be isotrivial (Proposition~\ref{prop:Liouville-general}). This criterion applies, in particular, when the base is an abelian variety (Remark~\ref{thm:LiouvillePeriods}).

When $V=X$ carries an algebraic action of a complex algebraic group $H$, each element $h\in H(\mathbb C)$ defines a pulled-back abelian variety $h^*G$. We compare its bad divisors with those of $G$ and use their invariance under $K$-isogeny to obtain rigidity under the action. The elliptic case illustrates the argument: for an elliptic curve over the function field of an elliptic base, isogenies preserve the places of bad reduction, whereas a non-torsion translation preserves no nonempty finite set of points. More generally, divisorial monodromy extends this obstruction: if $H$ has a dense open orbit in $X$ and a Zariski-dense subgroup of $H(\mathbb C)$ preserves the $K$-isogeny class of $G$, then $G$ is isotrivial (Corollary~\ref{thm:translation}). We also relate the groups of elements preserving the isomorphism or isogeny class of $G$ to the stabilizer of its bad divisors, describing their identity components and giving criteria for their finiteness.

\subsection{Divisorial bad reduction and isotriviality}\label{subsec:BadReduction}

Let $V$ be a smooth quasi-projective complex variety and put $K=\mathbb C(V)$. Let $G/K$ be an abelian variety of dimension $g$, equipped with a polarization. Choose a dense Zariski-open subset $U\subseteq V$ and an abelian scheme $\pi\colon\mathcal G\to U$ with generic fibre $G$. After shrinking $U$, we may assume that the chosen polarization extends to $\mathcal G$. The associated integral local system $\mathbb V:=R^1\pi_*\mathbb Z$ on $U(\mathbb C)$ has rank $2g$, with fibre $\mathbb V_q=H^1(\mathcal G_q,\mathbb Z)$ at $q\in U(\mathbb C)$. Write
\[
V\setminus U=D_1\cup\cdots\cup D_r\cup Y_2,
\]
where \(D_1,\ldots,D_r\) are the irreducible components of codimension one and \(Y_2\) is the union of the remaining components, of codimension at least two. We call the \(D_i\) the \emph{boundary divisors}. They depend on the choice of \(U\) and may include divisors at whose generic points \(G\) has good reduction.

Fix a base point \(q_0\in U(\mathbb C)\). Parallel transport along loops based at \(q_0\) defines the monodromy representation
\[
\rho:\pi_1(U(\mathbb C),q_0)
\longrightarrow \operatorname{GL}(\mathbb V_{q_0}).
\]
Let \(D\) be one of the boundary divisors \(D_1,\ldots,D_r\). Choose a point \(c\in D\) at which \(D\) is smooth and which belongs to no other irreducible component of \(V\setminus U\). Put \(d=\dim V\). In suitable holomorphic coordinates \((z_1,\ldots,z_d)\) centred at \(c\), a sufficiently small neighbourhood \(N\) is identified with a polydisc \(\Delta^d\), with
\[
(V\setminus U)\cap N= D\cap N=\{z_1=0\}, \qquad U\cap N=\Delta^*\times\Delta^{d-1}.
\]
For \(0<\varepsilon<1\), the positively oriented loop
\[
\gamma(t)=(\varepsilon e^{2\pi it},0,\ldots,0),
\qquad 0\le t\le1,
\]
generates \(\pi_1(U\cap N,(\varepsilon,0,\ldots,0))\). Choose a path \(\alpha\) in \(U\) from \(q_0\) to \((\varepsilon,0,\ldots,0)\). The homotopy class of the based loop
\[
\gamma_D=\alpha\cdot\gamma\cdot\alpha^{-1}
\]
is called a \emph{meridian} of $D$; we also denote this class by $\gamma_D$. Its conjugacy class in $\pi_1(U(\mathbb C),q_0)$ is independent of the point $c$, the local coordinates, and the path $\alpha$; see \cite[Definition~1.33 and Proposition~1.34]{rig:Cogolludo}. The \emph{local monodromy around \(D\)} is
\[
M_D:=\rho(\gamma_D)\in\operatorname{GL}(\mathbb V_{q_0}),
\]
defined up to conjugacy. After choosing a \(\mathbb Z\)-basis of \(\mathbb V_{q_0}\), we regard \(M_D\) as an element of \(\operatorname{GL}_{2g}(\mathbb Z)\). In particular, whether \(M_D\) is trivial, has finite order, or has infinite order is independent of these choices. For a prime divisor \(D\subseteq V\) which is not a boundary divisor, we put \(M_D=1\). This agrees with the same construction at a general point of \(D\): the transverse disc can then be chosen inside \(U\), so its boundary loop has trivial monodromy. Define
\[
\begin{aligned}
\operatorname{Bad}(G)
&:=\{D\subseteq V:\ D\text{ is a prime divisor and }M_D\ne1\},\\
\operatorname{Bad}^{\infty}(G)
&:=\{D\subseteq V:\ D\text{ is a prime divisor and }M_D
   \text{ has infinite order}\}.
\end{aligned}
\]
Both are finite sets of prime divisors of $V$, since only boundary divisors can occur. These sets are independent of the choice of \(U\) and of the abelian scheme \(\mathcal G\to U\). Indeed, shrinking \(U\) does not change the local monodromy around the original boundary divisors. Any additional boundary divisor has a general point in the original open subset, where the local system is defined on a whole neighbourhood; its local monodromy is therefore trivial. Since an identification of the generic fibres of two such abelian schemes extends to an isomorphism over a common dense open subset, the assertion follows. For $g=1$, $\operatorname{Bad}^{\infty}(G)$ is the set of prime divisors along which the $j$-invariant of $G$ has a pole.

\begin{rem}\label{rem:meridians-and-inertia}
The local monodromy around a divisor is closely related to the notion of bad reduction. Let $D\subseteq V$ be a prime divisor, with generic point $\eta_D$. Since $V$ is smooth, $R_D:=\mathcal O_{V,\eta_D}$ is a discrete valuation ring with fraction field $K$. Good reduction at $D$ means that $G$ extends to an abelian scheme over $\operatorname{Spec}R_D$; potentially good reduction means good reduction after a finite extension of $K$, at a valuation above $D$. We have
\[
\begin{aligned}
M_D=1 &\quad\Longleftrightarrow\quad G\text{ has good reduction at }D,\\
M_D\text{ has finite order} &\quad\Longleftrightarrow\quad G\text{ has potentially good reduction at }D.
\end{aligned}
\]

We briefly explain how these statements follow from the N\'eron--Ogg--Shafarevich criterion. Fix a prime $\ell$. The finite covers $\mathcal G[\ell^n]\to U$ describe the $\ell^n$-torsion points of the fibres. Write $N_D$ for the monodromy around $\gamma_D$ on $H_1(\mathcal G_{q_0},\mathbb Z)$. This action is dual to $M_D$, so the two automorphisms have the same order. For every $n\geq1$, the analytic description of an abelian variety gives
\[
\mathcal G_{q_0}[\ell^n] \simeq H_1(\mathcal G_{q_0},\mathbb Z)/ \ell^n H_1(\mathcal G_{q_0},\mathbb Z).
\]
Consequently, continuation around $\gamma_D$ acts on the $\ell^n$-torsion points by $N_D$ modulo $\ell^n$. This also describes the inertia action at $D$. Indeed, consider a finite Galois cover of $U$ over which these torsion points become single-valued, and extend it over $V$ by normalization. Near a general point of $D$, choose coordinates with $D=\{z_1=0\}$. Keeping the other coordinates fixed, the covering is locally $z_1=w_1^e$. Lifting a positive loop around $z_1=0$ sends $w_1$ to $e^{2\pi i/e}w_1$, giving a generator of the local inertia group. Thus, after compatible identifications, the inertia image on $G[\ell^n]$ is generated by $N_D$ modulo $\ell^n$. The Tate module records these actions for all $n$ simultaneously. Inertia acts trivially on it if and only if
\[
N_D\equiv I\pmod{\ell^n} \qquad\text{for every }n,
\]
which is equivalent to $N_D=I$. Similarly, its inertia image is finite if and only if there is a single integer $r\geq1$ such that
\[
N_D^r\equiv I\pmod{\ell^n} \qquad\text{for every }n,
\]
which is equivalent to $N_D^r=I$. Since $N_D$ and $M_D$ have the same order, this gives the required comparison. The N\'eron--Ogg--Shafarevich criterion \cite[Theorem~1]{S-T}, applied also after finite extensions, gives the asserted equivalences.
\end{rem}
The local system $\mathbb V$ can actually be extended across the good divisors. Put
\[
W:=V\setminus \bigcup_{D\in\operatorname{Bad}(G)}D.
\]
Thus $U\subseteq W$, and every boundary divisor whose generic point belongs to $W$ has trivial local monodromy. We show that this is precisely what is needed to extend the local system to $W$.

\begin{lem}\label{lem:extension-local-system}
The integral local system $\mathbb V$ extends to a local system $\mathbb V_W$ on $W(\mathbb C)$, uniquely up to a unique isomorphism compatible with its identification with $\mathbb V$ over $U$. In particular,
\[
\operatorname{Bad}(G)=\varnothing \quad\Longleftrightarrow\quad \mathbb V\text{ extends to a local system on }V(\mathbb C).
\]
\end{lem}

\begin{proof}
Recall that
\[
V\setminus U=D_1\cup\cdots\cup D_r\cup Y_2.
\]
Let
\[
Z:= Y_2 \cup\bigcup_{i=1}^r\operatorname{Sing}(D_i) \cup\bigcup_{1\leq i<j\leq r}(D_i\cap D_j), \qquad W^\circ:=W\setminus Z.
\]
The closed subset $Z$ has complex codimension at least two in $V$. Since $Z\subseteq V\setminus U$ and $U\subseteq W$, we have $U\subseteq W^\circ$. Moreover, $W^\circ\setminus U$ is a disjoint union of smooth divisors, all with trivial local monodromy by the definition of $W$. Fix $c\in W^\circ\setminus U$. There is a sufficiently small analytic polydisc $N\subseteq W^\circ$ centred at $c$ such that
\[
N\simeq\Delta^d, \qquad N\cap U\simeq\Delta^*\times\Delta^{d-1}, \qquad d=\dim V.
\]
The fundamental group of $N\cap U$ is generated by a meridian of the unique boundary divisor through $c$. Since its monodromy is trivial, $\mathbb V|_{N\cap U}$ is a constant local system and therefore extends to a constant local system on $N$. These extensions glue uniquely with $\mathbb V$. Indeed, a morphism between constant local systems on the connected punctured polydisc $N\cap U$ extends uniquely to $N$. Applying this observation locally on overlaps, the identifications over $U$ extend uniquely and satisfy the cocycle condition. We thus obtain a local system $\mathbb V^\circ$ on $W^\circ(\mathbb C)$.

It remains to extend across $Z\cap W$. We use the standard topological fact that removing a closed analytic subset of complex codimension at least two from a complex manifold does not change its fundamental group; the removed subset need not be smooth. Since $W$ is smooth, the inclusion $W^\circ\hookrightarrow W$ therefore induces an isomorphism
\[
\pi_1\bigl(W^\circ(\mathbb C),q_0\bigr) \xrightarrow{\;\sim\;} \pi_1(W(\mathbb C),q_0).
\]
The correspondence between local systems and representations of the fundamental group now gives an extension $\mathbb V_W$ to $W(\mathbb C)$, unique up to a unique isomorphism compatible with its restriction to $\mathbb V^\circ$. Together with the local uniqueness above, this proves the asserted uniqueness relative to $U$.

Finally, if $\operatorname{Bad}(G)=\varnothing$, then $W=V$, so the preceding construction gives an extension to $V(\mathbb C)$. Conversely, if $\mathbb V$ extends to $V(\mathbb C)$, every meridian has trivial monodromy, since its transverse circle bounds a disc in $V(\mathbb C)$. Hence $\operatorname{Bad}(G)=\varnothing$.
\end{proof}


The preceding lemma characterizes the absence of divisorial bad reduction in terms of extension of the integral local system. Under an additional hypothesis on the universal cover of the base, this extension forces isotriviality.

\begin{prop}\label{prop:Liouville-general}
Assume that every bounded holomorphic function on the analytic universal cover $\widehat V$ of $V(\mathbb C)$ is constant. If $\mathbb V=R^1\pi_*\mathbb Z$ extends to an integral local system $\widetilde{\mathbb V}$ on $V(\mathbb C)$, then $G$ is isotrivial. In particular, by Lemma~\ref{lem:extension-local-system}, this conclusion holds whenever $\operatorname{Bad}(G)=\varnothing$.
\end{prop}

\begin{proof}
Let $p\colon\widehat V\to V(\bC)$ be the universal covering and put $\Omega:=p^{-1}(U(\bC))$. The set $\Omega$ is connected, since it is the complement of a proper closed analytic subset of the connected complex manifold $\widehat V$. We use homology to describe the fibres as complex tori. The dual local system $\mathbb V^\vee$ has fibre $H_1(\mathcal G_u,\bZ)$ at $u\in U(\bC)$. Since $\widehat V$ is simply connected, $p^*(\widetilde{\mathbb V}^{\vee})$ is constant. Choose a flat trivialization, identifying its fibres with a fixed lattice $\Lambda$. The chosen polarization induces a flat nondegenerate alternating form on these homology lattices. Since $\Omega$ is connected, this form is represented in the chosen trivialization by a single alternating form $E \colon \Lambda \times \Lambda \longrightarrow \bZ$, whose extension to $\Lambda_{\bR}$ is nondegenerate.

The chosen flat trivialization identifies the complex structures of the fibres with complex structures on $\Lambda_{\mathbb R}$ compatible with $E$. After choosing a real symplectic basis for $E$, the space of such complex structures is identified with the Siegel upper half-space
\[
\mathfrak H_g := \left\{ \tau\in M_g(\mathbb C): \tau={}^t\tau,\ \operatorname{Im}\tau>0 \right\}.
\]
The holomorphic variation of the fibres therefore defines a holomorphic period map $\Psi \colon \Omega \longrightarrow \mathfrak H_g$. This map is single-valued because the pulled-back local system has been globally trivialized. The Cayley transform $\mathcal C(\tau) :=(\tau-iI_g)(\tau+iI_g)^{-1}$ identifies $\mathfrak H_g$ biholomorphically with a bounded domain. Consequently, every matrix entry of $\mathcal C\circ\Psi$ is a bounded holomorphic function on $\Omega$. Since $\widehat V\setminus\Omega$ is a proper closed analytic subset, the Riemann extension theorem \cite[Chapter~II, Remark~4.2]{rig:Demailly} extends these entries to holomorphic functions on $\widehat V$. Their extensions remain bounded by continuity and density, so they are constant by hypothesis. The injectivity of $\mathcal C$ then implies that $\Psi$ is constant.

Let $\lambda$ be the chosen polarization, of degree $\delta^2$, and let $\mu_\lambda\colon U\longrightarrow\bA_{g,\delta}$ be the associated classifying morphism to the coarse moduli space of polarized abelian varieties of dimension $g$ and polarization degree $\delta^2$. Since $\Psi$ is constant, all fibres $(\mathcal G_q,\lambda_q)$, for $q\in U(\bC)$, are isomorphic as polarized abelian varieties. Thus $\mu_\lambda$ sends every closed point of $U$ to the same moduli point. Since $U$ is reduced, $\mu_\lambda$ is constant as a morphism and so $\mathcal G/U$ is isotrivial.
\end{proof}

\begin{rem}\label{thm:LiouvillePeriods}
The hypothesis of Proposition~\ref{prop:Liouville-general} holds when $V$ is any connected complex algebraic group. Indeed, its analytic universal cover $\widehat V$ is a connected complex Lie group. If $f$ is a bounded holomorphic function on $\widehat V$, then, for every $a\in\widehat V$ and every vector $v$ in its Lie algebra, the function $z\mapsto f(a\exp(zv))$ is bounded and entire, hence constant. Differentiating at $z=0$ shows that $df$ vanishes in every tangent direction, so $f$ is constant.

The  case $\widehat V\simeq\mathbb C^d$ follows directly from Liouville's theorem and includes abelian and semiabelian varieties, as well as torsors under them. The hypothesis also holds for $V=A\times Y$, where $A$ is an abelian variety and $Y$ is a smooth simply connected projective variety: the universal cover is $\mathbb C^{\dim A}\times Y$, and a bounded holomorphic function on this product is constant along $Y$ by compactness and then along $\mathbb C^{\dim A}$ by Liouville's theorem.

For an abelian base, one can also specify a cover trivializing the family. Under the hypotheses of Proposition~\ref{prop:Liouville-general}, suppose that $V=A$ is an abelian variety. The proposition gives isotriviality, so the monodromy of $\mathbb V$ is finite. Since the inclusion $U(\bC)\hookrightarrow A(\bC)$ induces a surjection on fundamental groups, the monodromy group $\Gamma$ of $\widetilde{\mathbb V}$ is finite as well. Choose $n\geq1$ divisible by $|\Gamma|$ and put
\[
U_n:=[n]_A^{-1}(U), \qquad f_n:=[n]_A|_{U_n}\colon U_n\longrightarrow U.
\]
One can show that there is an isomorphism of polarized abelian schemes
\[
\mathcal G\times_{U,f_n}U_n \simeq B\times_{\bC}U_n
\]
for a polarized complex abelian variety $B$. The same conclusion holds for a torsor under an abelian variety, with multiplication defined using a chosen origin.
\end{rem}

\begin{rem}\label{rem:compact-base-isotriviality}
In relative dimension one, every abelian scheme over a smooth connected projective complex variety is isotrivial: its $j$-invariant defines a morphism to the affine line and is therefore constant. In higher relative dimension, compactness of the base alone no longer implies isotriviality. Indeed, non-isotrivial polarized abelian schemes over smooth projective curves exist already in relative dimension two, for example over compact quaternionic Shimura curves with sufficiently small level; see \cite[Theorem~1.1 and \S2.1]{rig:YuanKS}. Taking products with constant abelian varieties gives examples in every relative dimension $g\geq2$.

Taking $U=V$ in Proposition~\ref{prop:Liouville-general}, however, shows that every polarized abelian scheme over a smooth complex variety $V$ is isotrivial if the analytic universal cover of $V$ admits no nonconstant bounded holomorphic functions. This conclusion holds in every relative dimension and, in particular, when $V$ is a complex abelian variety.
\end{rem}

\subsection{Isogeny rigidity under algebraic group actions}
\label{subsec:translation-groups-K}

We retain the notation of \Cref{subsec:BadReduction} and write $X=V$, so that $K=\bC(X)$. Let $H$ be a complex algebraic group acting algebraically on $X$, with action morphism $H\times_{\bC}X\to X$ of $\bC$-schemes. The chosen polarization is auxiliary and will be used in the period-map and moduli arguments. Unless explicitly stated otherwise, homomorphisms, isomorphisms and isogenies are those of the underlying abelian varieties or abelian schemes, with no compatibility required with the chosen polarizations. For $h\in H(\mathbb C)$, we also write $h\colon X\to X$ for the automorphism $q\mapsto h\cdot q$. We denote by $h^*G$ the pullback of $G$ by the induced automorphism of $K$. Consider the abelian scheme
\[
h^*\mathcal G:=h^{-1}(U)\times_U\mathcal G \longrightarrow h^{-1}(U),
\]
where the fibre product is formed using $h|_{h^{-1}(U)}\colon h^{-1}(U)\to U$ and $\pi\colon\mathcal G\to U$. Its fibre at $q\in h^{-1}(U)(\mathbb C)$ is $\mathcal G_{h\cdot q}$, its generic fibre is $h^*G$, and its associated integral local system is $h^*\mathbb V$. To compare $\mathcal G$ and $h^*\mathcal G$ and their associated local systems, we restrict both schemes to the common dense open subset $U\cap h^{-1}(U)$. This does not change their generic fibres or the associated bad-divisor sets. Define the following subgroups of $H(\mathbb C)$:
\[
\begin{aligned}
\operatorname{Isom}_H(G) &:=\{h\in H(\mathbb C):h^*G\simeq_K G\},\\
\operatorname{Isog}_H(G) &:=\{h\in H(\mathbb C): h^*G\text{ is }K\text{-isogenous to }G\}.
\end{aligned}
\]
Notice that $\operatorname{Isom}_H(G)\subseteq\operatorname{Isog}_H(G)$. These definitions concern the underlying abelian variety $G$ and are independent of the auxiliary polarization. We call these sets the \emph{isomorphism group} and the \emph{isogeny group} of $G$ relative to the $H$-action, respectively. When $X=A$ is an abelian variety and $H=A$ acts by translations $T_x\colon q\mapsto q+x$, we call such groups the \emph{translation--isomorphism group} and the \emph{translation--isogeny group}, denoted by $\operatorname{Isom}(G)$ and $\operatorname{Isog}(G)$.

We first show that the bad-divisor sets are invariant under $K$-isogenies. Let $\pi'\colon\mathcal G'\to U$ be another abelian scheme with generic fibre $G'$, and let $\phi\colon G\to G'$ be a $K$-isogeny. Since $U$ is normal, $\phi$ extends uniquely to an isogeny of abelian schemes $\Phi\colon\mathcal G\longrightarrow\mathcal G'$ over $U$. Put $\mathbb V':=R^1\pi'_*\mathbb Z$. The induced morphism on cohomology gives an isomorphism of rational local systems
\[
\Phi_{\mathbb Q}^{*}\colon \mathbb V'\otimes\mathbb Q \xrightarrow{\;\sim\;} \mathbb V\otimes\mathbb Q.
\]
Write $M_D(G)$ and $M_D(G')$ for the local monodromies of $G$ and $G'$ around a prime divisor $D$, respectively. Since $\Phi_{\mathbb Q}^{*}$ commutes with parallel transport, these monodromies are conjugate in $\operatorname{GL}_{2g}(\mathbb Q)$ for every prime divisor $D\subseteq X$. Conjugacy preserves triviality and finiteness of order, so
\begin{equation}\label{eq:bad-isogeny-invariance}
    \operatorname{Bad}(G)=\operatorname{Bad}(G'), \qquad \operatorname{Bad}^{\infty}(G) =\operatorname{Bad}^{\infty}(G').
\end{equation}

We next consider how the bad-divisor sets behave under the action on the base. The group $H(\bC)$ acts on the power set of the prime divisors of $X$, i.e. for $h\in H(\bC)$ and a set $\mathfrak D$ of prime divisors, we have $h\cdot\mathfrak D:=\{h(D):D\in\mathfrak D\}$. Since $h$ sends a meridian of $h^{-1}(D)$ to a meridian of $D$, the local monodromies $M_{h^{-1}(D)}(h^*G)$ and $M_D(G)$ are conjugate. Therefore, we obtain
\begin{equation}\label{eq:bad-translation}
    \operatorname{Bad}(h^*G)=h^{-1}\cdot\operatorname{Bad}(G), \qquad \operatorname{Bad}^{\infty}(h^*G)=h^{-1}\cdot\operatorname{Bad}^{\infty}(G).
\end{equation}
We now consider the stabilizer $\Sigma(\mathfrak D)$ of a set $\mathfrak D$ of prime divisors with respect to the previous action. In particular, we write
\[
\Sigma_G:=\Sigma(\operatorname{Bad}(G)), \qquad \Sigma_G^\infty:=\Sigma(\operatorname{Bad}^{\infty}(G)), \qquad \Sigma(D):=\Sigma(\{D\}), \qquad \Sigma(\varnothing):=H.
\]
\begin{lem}\label{lem:divisor-stabilizers}
For any finite set $\mathfrak D$ of prime divisors of $X$, the group $\Sigma(\mathfrak D)$ is a closed algebraic subgroup of $H$, and its identity component $\Sigma(\mathfrak D)^0$ is a connected closed algebraic subgroup satisfying
\[
\Sigma(\mathfrak D)^0 = \left( \bigcap_{D\in\mathfrak D}\Sigma(D) \right)^0.
\]
\end{lem}

\begin{proof}
The assertion is immediate if $\mathfrak D$ is empty, with the empty intersection interpreted as $H$. Otherwise, write $\mathfrak D=\{D_1,\ldots,D_r\}$ and put $Y:=D_1\cup\cdots\cup D_r$. The automorphism $h$ gives a bijection between the irreducible components of $Y$ and $h(Y)$. Since these components are prime divisors, $h\in\Sigma(\mathfrak D)$ is equivalent to the containment $h(Y)\subseteq Y$, and we have
\[
\Sigma(\mathfrak D)(\mathbb C)
=\{h\in H(\mathbb C):h(Y)\subseteq Y\}
=\bigcap_{p\in Y(\mathbb C)}a_p^{-1}(Y)(\mathbb C),
\]
where $a_p\colon H\to X$, $h\mapsto h\cdot p$, is the orbit map. Therefore, $\Sigma(\mathfrak D)$ is a closed algebraic subgroup of $H$, considered with its reduced structure. The subgroup $L:=\bigcap_{i=1}^r\Sigma(D_i)$ is closed and is the kernel of the permutation action on $\mathfrak D$. It therefore has finite index in $\Sigma(\mathfrak D)$. Its cosets are closed and, since there are only finitely many, also open in $\Sigma(\mathfrak D)$. Consequently $L$ and $\Sigma(\mathfrak D)$ have the same identity component.
\end{proof}

\begin{prop}\label{inv_isog}
$\operatorname{Isom}_H(G) \subseteq\operatorname{Isog}_H(G) \subseteq \Sigma_G(\mathbb C) \cap \Sigma_G^\infty(\mathbb C)$.
\end{prop}

\begin{proof}
Combine \eqref{eq:bad-isogeny-invariance} and \eqref{eq:bad-translation}.
\end{proof}

We denote by $\Sigma_G^0$ the identity component of $\Sigma_G$. When no confusion arises, we use the same notation for an algebraic subgroup of $H$ and its group of complex points. We first prove isotriviality along the Zariski closure of the $\Sigma_G^0$-orbits.

\begin{lem}\label{lem:coset-isotriviality}
For every $\Sigma_G^0$-orbit $\Xi$ in $X$ such that $\Xi\cap U\ne\varnothing$, the restriction $\mathcal G|_{\overline{\Xi}^{\mathrm{Zar}}\cap U}$ is isotrivial. 
\end{lem}

\begin{proof}
Recall the open subset
\[
W=X\setminus \bigcup_{D\in\operatorname{Bad}(G)}D
\]
defined by the bad divisors. By Lemma \ref{lem:divisor-stabilizers} every bad divisor is stable under $\Sigma_G^0$, so the subset $W$ is a union of $\Sigma_G^0$-orbits. Since $U\subseteq W$, this implies $\Xi\subseteq W$. Let $\mathbb V_W$ be the extension of $\mathbb V$ to $W(\bC)$ provided by Lemma~\ref{lem:extension-local-system}. Choose $q\in\Xi\cap U$ and consider the surjective orbit map
\[
a_q\colon\Sigma_G^0\longrightarrow\Xi,\qquad h\longmapsto h\cdot q.
\]
The local system of $a_q^*\mathcal G$ on $a_q^{-1}(U)$ extends over $\Sigma_G^0$ as $a_q^*\mathbb V_W$. Proposition~\ref{prop:Liouville-general}, applied as in Remark~\ref{thm:LiouvillePeriods}, shows that the pulled-back family is isotrivial. Since $a_q$ is surjective, the polarized moduli map of $\mathcal G$ is constant on $\Xi\cap U$. This subset is dense in $\overline{\Xi}^{\mathrm{Zar}}\cap U$, so the moduli map is constant on $\overline{\Xi}^{\mathrm{Zar}}\cap U$ as well. Thus the restricted family is isotrivial.
\end{proof}

We have seen that every element of $H(\mathbb C)$ preserving the $K$-isogeny class of $G$ belongs to $\Sigma_G$. We now show that every point of $\Sigma_G^0$ induces an automorphism of the base preserving the $K$-isomorphism class of $G$. More precisely, the corresponding isomorphisms can be chosen to preserve the auxiliary polarization $\lambda$.

\begin{thm}\label{rig:thm:invariance}
Let us denote by $\lambda$ a fixed polarization of the abelian scheme $\mathcal G/U$. For every $h\in \Sigma_G^0$, there is an isomorphism of polarized abelian schemes
\[
(h^*\mathcal G,h^*\lambda)|_{U_h} \simeq (\mathcal G,\lambda)|_{U_h}, \qquad U_h:=U\cap h^{-1}(U).
\]
In particular at the level of complex points we have
\begin{equation}\label{eq:translation-group-structure}
    \Sigma_G^0\subseteq\operatorname{Isom}_H(G) \subseteq\operatorname{Isog}_H(G) \subseteq\Sigma_G,
\end{equation}
and both $\operatorname{Isom}_H(G)$ and $\operatorname{Isog}_H(G)$ are  closed algebraic subgroups of $H$, with common identity component $\Sigma_G^0$.
\end{thm}

\begin{proof}
Let $\Xi$ be a $\Sigma_G^0$-orbit intersecting $U$;  Lemma~\ref{lem:coset-isotriviality} shows that $\mathcal G|_{\Xi\cap U}$ is isotrivial. Its Hodge filtration is therefore locally constant in flat trivializations. Since the polarization form is also flat, parallel transport along paths in $\Xi\cap U$ preserves the integral Hodge structures and the polarization forms. Notice that this conclusion remains true along any path in $\Xi(\mathbb C)$ with endpoints in $\Xi\cap U$. Indeed, since $\Xi\setminus U$ is a proper closed analytic subset of the complex manifold $\Xi$, any such path  is homotopic  to a path contained in $\Xi\cap U$, and parallel transport in $\mathbb V_W|_{\Xi(\bC)}$ depends only on the homotopy class of the path.

Fix $h\in\Sigma_G^0$ and choose a continuous path
\[
\eta\colon[0,1]\longrightarrow\Sigma_G^0,
\qquad
\eta(0)=e,\quad \eta(1)=h,
\]
where $e$ denotes the identity element of the group. Since $W$ is a union of $\Sigma_G^0$-orbits (see the proof of Lemma~\ref{lem:coset-isotriviality}), it is stable under $\Sigma_G^0$. So the map
\[
H_\eta\colon W(\bC)\times[0,1]\longrightarrow W(\bC),
\qquad
H_\eta(q,t)=\eta(t)\cdot q,
\]
is a homotopy from the identity map to $h|_W$. Parallel transport along the paths $t\mapsto \eta(t)\cdot q$ defines an isomorphism of integral local systems
\[
\theta_\eta\colon
\mathbb V_W\xrightarrow{\;\sim\;}h^*\mathbb V_W.
\]
Indeed, the homotopy $H_\eta$ ensures that these fibrewise maps commute with parallel transport along paths in $W(\bC)$. For $q\in U_h$, the path $t\mapsto \eta(t)\cdot q$ lies in the orbit $\Sigma_G^0\cdot q$ and joins the two points $q,h\cdot q\in U$. By the preceding argument, the fibre map
\[
(\theta_\eta)_q\colon
H^1(\mathcal G_q,\bZ)
\xrightarrow{\;\sim\;}
H^1(\mathcal G_{h\cdot q},\bZ)
\]
preserves the Hodge structures and the polarization forms. Thus $\theta_\eta|_{U_h}$ is an isomorphism of the polarized integral variations of Hodge structures associated with $\mathcal G|_{U_h}$ and $h^*\mathcal G|_{U_h}$. By the full faithfulness of the contravariant functor associating to an abelian scheme its integral variation of Hodge structures \cite[(4.4.3)]{rig:Deligne}, $\theta_\eta|_{U_h}$ is induced by a unique algebraic isomorphism
\[
f_\eta\colon
h^*\mathcal G|_{U_h}
\xrightarrow{\;\sim\;}
\mathcal G|_{U_h}
\]
compatible with the polarization.
This gives the inclusion $\Sigma_G^0\subseteq\operatorname{Isom}_H(G)$. Together with  Proposition~\ref{inv_isog}, this proves \eqref{eq:translation-group-structure}.

Finally, $\Sigma_G$ has finitely many connected components, each a coset of $\Sigma_G^0$. Every subgroup lying between $\Sigma_G^0$ and $\Sigma_G$ is a union of these cosets so the last assertion follows immediately.
\end{proof}

We now deduce the isogeny rigidity result. Its translation case provides the key rigidity argument for the proof of Theorem~\ref{thm:B}.

\begin{cor}\label{thm:translation}
If $\operatorname{Isog}_H(G)$ has a Zariski-dense orbit in $X$, then $G$ is isotrivial. In particular, if $h\in H(\mathbb C)$ is such that $\overline{\langle h\rangle}^{\mathrm{Zar}}$ has a Zariski-dense orbit in $X$ and $h^k\in \operatorname{Isog}_H(G)$ for some $k\in\mathbb Z\setminus\{0\}$, then $G$ is isotrivial.
\end{cor}

\begin{proof}
Put $L:=\operatorname{Isog}_H(G)$. For every $x\in U(\mathbb C)$, put
\[
Y_x:=\overline{L^0\cdot x}^{\mathrm{Zar}}.
\]
Then $\mathcal G|_{Y_x\cap U}$ is isotrivial by Lemma ~\ref{lem:coset-isotriviality} and  Theorem  ~\ref{rig:thm:invariance}. Suppose that $L$ has a Zariski-dense orbit $O$ in $X$. Orbits of algebraic group actions are locally closed, hence open in their closures. Therefore $O$ is open in $X$. Since $L/L^0$ is finite, the $L^0$-orbits in $O$ have the same dimension as $O$ and are also open in $O$. Two nonempty open subsets of the irreducible variety $O$ cannot be disjoint, so $L^0$ acts transitively on $O$. Choosing $x\in O\cap U$, we obtain $Y_x=X$.

For the assertion concerning a single element, put
\[
M:=\overline{\langle h\rangle}^{\mathrm{Zar}},
\qquad
M_k:=\overline{\langle h^k\rangle}^{\mathrm{Zar}}.
\]
The group $M$ is a union of finitely many cosets of $M_k$, so $M_k^0=M^0$ and $M_k$ also has a dense orbit in $X$. Since $h^k\in L$ and $L$ is closed, we have $M_k\subseteq L$. Thus $L$ has a dense orbit in $X$, and the preceding conclusion proves that $G$ is isotrivial.
\end{proof}

\begin{rem}\label{rem:group-action-examples}
The following examples illustrate the global conclusion of Corollary~\ref{thm:translation}. In each case, $G$ is an abelian variety over $K=\mathbb C(X)$, and all isogenies are defined over $K$.
\begin{enumerate}
\item[$(i)$] Let $X=H=A$ be an abelian variety acting on itself by translations. If
\[
\overline{\mathbb Zx}^{\mathrm{Zar}}=A
\qquad\text{and}\qquad
kx\in\operatorname{Isog}(G)
\quad\text{for some }k\in\mathbb Z\setminus\{0\},
\]
then $G$ is isotrivial, since $\mathbb Z(kx)$ is still Zariski dense in $A$. This is the criterion used below in the proof of Theorem~\ref{thm:B}. If $A$ is simple of positive dimension, the density condition is equivalent to $x$ being non-torsion. The same conclusion holds for a torsor $X$ under $H=A$.

\item[$(ii)$] Let $X=H=E_1\times E_2$, acting by translations, where $E_1,E_2$ are elliptic curves. Choose non-torsion points $p_i\in E_i(\mathbb C)$. If
\[
(p_1,0),\ (0,p_2)\in\operatorname{Isog}(G),
\]
then $G$ is isotrivial, since these points generate the Zariski-dense subgroup $\mathbb Zp_1\times\mathbb Zp_2$ of $H$. One condition alone does not suffice: a non-isotrivial family pulled back from $E_2$ is preserved by all translations in the $E_1$-factor.

\item[$(iii)$] Let $X=\mathbb A^m$ and $H=\mathbb G_a^m$, acting by translations. For $i=1,\ldots,m$, let
\[
T_i(t_1,\ldots,t_m)=(t_1,\ldots,t_i+1,\ldots,t_m).
\]
If $T_i^*G$ is isogenous to $G$ for every $i$, then $G$ is isotrivial. Indeed, the corresponding elements of $H(\mathbb C)$ generate the Zariski-dense subgroup $\mathbb Z^m$, and $H$ acts transitively on $X$. For $m>1$, a single simultaneous translation by $(1,\ldots,1)$ does not suffice: its cyclic subgroup has closure equal to the diagonal copy of $\mathbb G_a$.

\item[$(iv)$] Let $X=H=\mathbb G_m^m$, acting by multiplication, and choose multiplicatively independent numbers $q_1,\ldots,q_m\in\mathbb C^*$. Put
\[
h(t_1,\ldots,t_m)=(q_1t_1,\ldots,q_mt_m).
\]
If $h^*G$ is isogenous to $G$, then $G$ is isotrivial, since the element $(q_1,\ldots,q_m)$ generates a Zariski-dense subgroup of $H$. For $m=1$, the condition is simply that $q_1$ is not a root of unity.

\item[$(v)$] Let $X=E\times\mathbb P^1$ and $H=E\times\mathbb G_m$, where $E$ is an elliptic curve, with action
\[
(b,s)\cdot(a,t)=(a+b,st).
\]
Choose a non-torsion point $p\in E(\mathbb C)$ and $q\in\mathbb C^*$ which is not a root of unity. If both
\[
h_1(a,t)=(a+p,t),\qquad h_2(a,t)=(a,qt)
\]
belong to $\operatorname{Isog}_H(G)$, then $G$ is isotrivial. Indeed, the corresponding elements $(p,1)$ and $(0,q)$ generate a Zariski-dense subgroup of $H$, whose dense open orbit is $E\times\mathbb G_m$.

\item[$(vi)$] Let $X=\mathbb A^2$ and $H=\operatorname{SL}_2$, acting in the natural way. If both
\[
h_1(u,v)=(u+v,v),\qquad h_2(u,v)=(u,v+u)
\]
belong to $\operatorname{Isog}_H(G)$, then $G$ is isotrivial. Indeed, the Zariski closure of the subgroup they generate contains the upper and lower unipotent subgroups, hence all of $\operatorname{SL}_2$. The latter has the dense open orbit $\mathbb A^2\setminus\{0\}$.
\end{enumerate}
\end{rem}

Put $d:=\dim X$, and let $\mu_\lambda$ be the coarse
moduli map of the polarized abelian scheme
$(\mathcal G,\lambda)\to U$. Define
\[
  \operatorname{var}(G)
  :=\dim\overline{\mu_\lambda(U)}^{\mathrm{Zar}}.
\]
This number is independent of the chosen model and
polarization. We say that $G$ has \emph{maximal variation} if
$\operatorname{var}(G)=d$, equivalently if
$\mu_\lambda$ is generically finite onto its image.

\begin{prop}\label{prop:translation-group-finiteness}
For every $x\in X(\mathbb C)$, we have
\[
\dim(\Sigma_G^0\cdot x)+\operatorname{var}(G) \leq d.
\]
In particular, if the $H$-action has finite kernel and $G$ has maximal variation, then $\operatorname{Isom}_H(G)$ and $\operatorname{Isog}_H(G)$ are finite.
\end{prop}

\begin{proof}
Let $r$ be the dimension of a general $\Sigma_G^0$-orbit, which is also the maximal orbit dimension. By Lemma~\ref{lem:coset-isotriviality}, $\mu_\lambda$ is constant on the intersection of each  $\Sigma_G^0$-orbit with $U$. Its general fibres therefore have dimension at least $r$, giving
\[
r+\operatorname{var}(G)\leq d.
\]

If $G$ has maximal variation, then $r=0$. Since $\Sigma_G^0$ is connected, all its orbits are singletons, so $\Sigma_G^0$ lies in the kernel of the action. A finite kernel therefore forces $\Sigma_G^0=\{e\}$, and the conclusion follows from Theorem~\ref{rig:thm:invariance}.
\end{proof}

We now specialize to $X=H=A$, where $A$ is an abelian variety of dimension $d\geq1$ acting by translations. Let
\[
q\colon A\longrightarrow Q:=A/\Sigma_G^0, \qquad d':=\dim Q.
\]
By Theorem~\ref{rig:thm:invariance}, $\operatorname{Isom}(G)$, $\operatorname{Isog}(G)$ and $\Sigma_G$ have common identity component $\Sigma_G^0$. Thus their quotients by $\Sigma_G^0$ are finite groups whose orders count their connected components.

By Lemma~\ref{lem:divisor-stabilizers}, every $D\in\operatorname{Bad}(G)$ is invariant under $\Sigma_G^0$, hence is a union of fibres of $q$. Since $q$ is proper and smooth with connected fibres, its image $\overline D:=q(D)$, endowed with the reduced structure, is a prime divisor on $Q$ and $D=q^*\overline D$. Distinct bad divisors have distinct images. In particular, if $\operatorname{Bad}(G)\ne\varnothing$, then $d'\geq1$.

\begin{prop}\label{prop:translation-group-bound}
Suppose that $\operatorname{Bad}(G)\ne\varnothing$,
and define the reduced divisors
\[
  \Delta_G:=\sum_{D\in\operatorname{Bad}(G)}D,
  \qquad
  \overline\Delta_G
  :=\sum_{D\in\operatorname{Bad}(G)}\overline D.
\]
Then $\Delta_G=q^*\overline\Delta_G$, and
$\overline\Delta_G$ is ample on $Q$. The orders of
$\operatorname{Isom}(G)/\Sigma_G^0$ and
$\operatorname{Isog}(G)/\Sigma_G^0$ divide
the positive integer $\frac{\overline\Delta_G^{\,d'}}{d'!}$. If at least one $\overline D$ is ample on $Q$,
then also
\[
\bigl|\operatorname{Isog}(G)/\Sigma_G^0\bigr| \leq |\operatorname{Bad}(G)| \min_{\substack{D\in\operatorname{Bad}(G)\\ \overline D\ \mathrm{ample}}} \frac{\overline D^{\,d'}}{d'!}.
\]
Moreover, $\operatorname{Isom}(G)$ and
$\operatorname{Isog}(G)$ are finite if and only if
$\Delta_G$ is ample.
\end{prop}

\begin{proof}
The equality $\Delta_G=q^*\overline\Delta_G$ follows
by summing the pullback identities for the individual
bad divisors. Since $\Sigma_G$ is the translation
stabilizer of $\Delta_G$, the translation stabilizer
of $\overline\Delta_G$ on $Q$ is the finite group $F:=\Sigma_G/\Sigma_G^0$. An effective divisor on an abelian variety is ample if and only if its translation stabilizer is finite; see \cite[\S6, Application~1]{Mum}. Thus $\overline\Delta_G$ is ample.

Put $E:=\overline\Delta_G$, and consider the finite
\'etale quotient $\rho\colon Q\longrightarrow Q/F$. Endow $\rho(E)$ with its reduced structure.
Since $E$ is reduced and $F$-invariant, we have
$E=\rho^*\rho(E)$. The projection formula and
Riemann--Roch \cite[\S16]{Mum} give
\[
  \frac{E^{d'}}{d'!}
  =
  |F|\frac{\rho(E)^{\,d'}}{d'!}
  =
  |F|\,\chi\bigl(\mathcal O_{Q/F}(\rho(E))\bigr).
\]
The left-hand side is positive because $E$ is ample,
so the Euler characteristic on the right is a
positive integer. Hence $|F|$ divides $E^{d'}/d'!$
and is at most this number. The inclusions in
Theorem~\ref{rig:thm:invariance} induce
\[
  \operatorname{Isom}(G)/\Sigma_G^0
  \subseteq
  \operatorname{Isog}(G)/\Sigma_G^0
  \subseteq F.
\]
This proves the divisibility assertions.

Now let $D\in\operatorname{Bad}(G)$ be such that
$\overline D$ is ample, and denote its finite
translation stabilizer on $Q$ by $F'$.
The group $F$ permutes the prime components of
$\overline\Delta_G$, so the orbit of $\overline D$
has at most $|\operatorname{Bad}(G)|$ elements.
Its stabilizer in $F$ is $F\cap F'$, and
orbit--stabilizer therefore gives
\[
\bigl|\operatorname{Isog}(G)/\Sigma_G^0\bigr| \leq |F| = |F\cdot\overline D|\,|F\cap F'| \leq |\operatorname{Bad}(G)|\,|F'|.
\]
Applying the same quotient and Riemann--Roch
argument to $\overline D$ and $F'$ yields
\[
  |F'|\leq\frac{\overline D^{\,d'}}{d'!}.
\]
Taking the minimum over the ample $\overline D$
proves the claim.

Finally, since $\operatorname{Isom}(G)$,
$\operatorname{Isog}(G)$ and $\Sigma_G$ have the
same identity component, the first two groups
are finite if and only if $\Sigma_G$ is finite.
As $\Sigma_G$ is the translation stabilizer of
$\Delta_G$, the ampleness criterion used above
shows that this is equivalent to $\Delta_G$
being ample.
\end{proof}

%% file: DoubleFamilies.tex
\subsection{Motivation: finite orbit points on a family}

\begin{defn}\label{defn:singleFamily}
Let $g\geq 1$. A \emph{family of $g$-dimensional abelian varieties} consists of a dominant rational map $f:\mathcal{A}'\dashrightarrow S'$ from a nonsingular projective variety to a  projective variety, satisfying the following conditions. Writing $D_f:=\mathcal{A}'\setminus\operatorname{Fund}(f)$, there exists a dense open subset $U\subseteq S'$ such that:
\begin{itemize}
    \item[$(i)$] the restriction $(f|_{D_f})^{-1}(U)\longrightarrow U$ is a smooth projective surjective morphism,

    \item[$(ii)$] all the fibres over $U$ are isomorphic, as varieties, to abelian varieties of dimension $g$.
\end{itemize}
We denote by $S$ the largest open subset of $S'$ satisfying these conditions, and set $\mathcal{A}:=(f|_{D_f})^{-1}(S)$. We call $S$ the \emph{smooth base} and $f:\mathcal{A}\to S$ the \emph{smooth part} of the family.
\end{defn}

The variety $S$ is nonsingular as well. Indeed, $\cA$ is nonsingular and  smoothness over $\mathbb{C}$ descends through the smooth surjective morphism $f:\mathcal{A}\to S$; see \cite[Tag~02K5]{stacks-project}. Writing $d:=\dim S'=\dim S$, we have
\[
\dim\mathcal{A}'=\dim\mathcal{A}=g+d.
\]

\begin{rem}
\label{rem:properness}
Properness of $f\colon\cA\to S$ is included in condition $(i)$ of Definition~\ref{defn:singleFamily}. In particular, for every $s\in S$, the fibre $\cA_s$ is proper over $\bC$. Since $\cA'$ is separated, its inclusion into $\cA'$ has closed image; see \cite[Tag~01W6]{stacks-project}. Thus $\cA_s$ is already closed in $\cA'$, and its closure does not acquire points of $\Fund(f)$. This requirement is not automatic for a rational map, even when the closures of its general fibres are abelian varieties. For example, contracting a $(-1)$-section of a rational elliptic
surface gives a rational map whose general fibre on its domain of definition is a punctured elliptic curve. Shrinking the base does not make that restriction proper.

When checking condition $(i)$, however, it suffices to establish properness over some dense open subset of $S'$. Indeed, suppose that $U\subseteq S'$ is dense open and
\[
f_U\colon(f|_{D_f})^{-1}(U)\longrightarrow U
\]
is proper. Since $f_U$ is dominant and proper, it is surjective. Let $N_U$ be its nonsmooth locus. Then
\[
U^{\mathrm{sm}}:=U\setminus f_U(N_U)
\]
is the largest open subset of $U$ over which $f_U$ is smooth. It is open by properness, and dense by generic smoothness in characteristic zero, since the source of $f_U$ is nonsingular. The restriction over $U^{\mathrm{sm}}$ is therefore smooth and surjective. It is also projective, since it is proper and its source is quasi-projective. Thus condition $(i)$ holds after replacing $U$ by $U^{\mathrm{sm}}$.
\end{rem}

Let $f:\cA\to S$ be the smooth part of a family of $g$-dimensional abelian varieties. Although each fibre is isomorphic to an abelian variety, no origins on the fibres have been chosen, and $f$ need not admit a global section. Thus \(\mathcal A\to S\) is not assumed to carry an abelian-scheme structure. Nevertheless, the family has a canonical abelian scheme of translations. Indeed, since \(f\) is projective and smooth, the relative automorphism functor
\[
T \mapsto \Aut_T(\cA \times_S T)
\]
is represented by a separated $S$-group scheme $\underline{\Aut}_{\cA/S}$; see
\cite[Exp.~221, \S4(c)]{FGA}. Its relative identity component
\[
\cB:=\underline{\Aut}_{\cA/S}^0
\]
is an abelian scheme over $S$ of relative dimension $g$. For each $s\in S$, its fibre $\cB_s$ is the identity component of the automorphism scheme of $\cA_s$ and consists precisely of the translations of $\cA_s$. We call $\cB/S$ the \emph{associated abelian scheme}, or the \emph{abelian scheme of translations}, of $\cA\to S$. The local justification of these assertions is given in Remark~\ref{rem:localOrigins} below.

Evaluation defines a canonical action, which we write additively:
\[
\mu:\cB\times_S\cA \longrightarrow\cA, \qquad (b,p)\longmapsto p + b.
\]
This action is simply transitive, in the sense that the morphism
\[
\theta:\cA \times_S \cB \longrightarrow \cA \times_S \cA, \qquad (p,b)\longmapsto(p,p + b),
\]
is an isomorphism. Its inverse has the form $\theta^{-1}(p,q)=(p,\delta(p,q))$, where $\delta:\mathcal{A}\times_S\mathcal{A}\to\mathcal{B}$ is the \emph{difference morphism}, characterized by
$p + \delta(p,q) = q$. For $p,q\in\cA_s(\bC)$, we therefore write $q-p:=\delta(p,q)\in\cB_s$. Thus two points of the same fibre have a canonical difference
in $\cB_s$, even though $\cA_s$ has no distinguished origin.

\begin{rem}\label{rem:localOrigins}
    Since $f$ is smooth and surjective, there exists an \'etale covering $\{U_i\to S\}$ by connected schemes such that each pullback $\cA_{U_i}\to U_i$ admits a section $o_i$; see \cite[Tag~055U]{stacks-project}. By \cite[Theorem~6.14]{MFK}, each $\cA_{U_i}/U_i$ has a unique abelian-scheme structure with identity section $o_i$. For this structure, the relative identity component of the
    automorphism scheme consists of translations. Hence the orbit map
    \[
    \cB_{U_i}\xrightarrow{\;\sim\;}\cA_{U_i}, \qquad b\longmapsto o_i+b,
    \]
    is an isomorphism identifying $\mu$ with addition on $\cB_{U_i}$. This proves, by \'etale descent, that $\cB/S$ is an abelian scheme of relative dimension $g$ and that $\theta$ is an isomorphism, and exhibits $\cA$ as a \emph{torsor under $\cB$ for the \'etale topology}\footnote{Since $\cB/S$ is smooth, the notions of \'etale torsor and fppf torsor under $\cB$ coincide.}. In particular, a global section $o\colon S\to\cA$ identifies $\cA$ with $\cB$, giving $\cA/S$ an abelian-scheme structure with identity section $o$. These identifications depend on the chosen origins, whereas $\cB$, its action, and the difference morphism are canonical. More precisely, any other abelian scheme over $S$ acting on
    $\cA$ and making it a torsor admits a unique $S$-group scheme isomorphism to $\cB$ compatible with the two actions.
\end{rem}

\begin{defn}\label{def:translation}
    A \emph{translation of $\cA/S$} is an $S$-automorphism $t: \cA\to\cA$ such that, for every $s\in S$, there exists $b_s\in\cB_s$ satisfying
    \[
    t(p)=p+b_s \qquad \textrm{for every } p\in\cA_s.
    \]
    The point $b_s$ is necessarily unique. We denote the group of translations of $\cA/S$ by
    $\operatorname{Trans}_S(\cA)$.
\end{defn}

By the representability of the relative automorphism functor, an $S$-automorphism $t$ of $\cA$ corresponds to a unique section
\[
\widetilde t\colon S\longrightarrow \underline{\Aut}_{\cA/S},
\]
whose value at $s\in S$ represents the restriction of $t$ to $\cA_s$. The subgroup scheme
$\cB\subseteq\underline{\Aut}_{\cA/S}$ is closed: its inclusion is a proper monomorphism, since $\cB/S$ is proper and $\underline{\Aut}_{\cA/S}$ is separated over $S$. As $S$ is reduced and its closed points are dense, the fibrewise condition in Definition~\ref{def:translation} is therefore equivalent to $\widetilde t$ factoring through
$\cB$. Consequently, evaluation induces a canonical isomorphism of groups
\[
\cB(S) \xrightarrow{\;\sim\;} \operatorname{Trans}_S(\cA), \qquad \sigma\longmapsto t_\sigma,
\]
where
\[
t_\sigma\colon\cA\longrightarrow\cA, \qquad p\longmapsto p+\sigma(f(p)).
\]
We call the section $\sigma$ the \emph{translation vector} of $t_\sigma$. Concretely, the translation vector of $t$ is characterized by
\[
\sigma(f(p))=\delta(p,t(p))=t(p)-p \qquad \textrm{for every }p\in\cA.
\]

Since $\cA$ is a dense open subset of $\cA'$, every translation $t_\sigma\colon\cA\to\cA$ determines a unique birational self-map $t_\sigma\colon\cA'\dashrightarrow\cA'$, which we denote by the same symbol, whose inverse is induced by $t_{-\sigma}$. Consequently, the identification of sections with translations gives an injective group homomorphism
\begin{equation}\label{eq:translationsBirational}
\cB(S)\hookrightarrow\operatorname{Bir}(\cA'), \qquad \sigma\longmapsto t_\sigma,
\end{equation}
where $\operatorname{Bir}(\cA')$ denotes the group of birational self-maps of $\cA'$.

\begin{rem}\label{rem:equivalentFamily}

    Let $g\geq1$, and let $f\colon\cA'\dashrightarrow S'$ be a dominant rational map from a nonsingular projective variety to a projective variety. Put $D_f:=\cA'\setminus\Fund(f)$, and assume that condition $(i)$ of Definition~\ref{defn:singleFamily} holds over some dense open subset of $S'$. Let $S\subseteq S'$ be the largest open subset over which $f\colon\cA:=(f|_{D_f})^{-1}(S)\longrightarrow S$ is smooth, projective, and surjective. Then $S$ is nonsingular. Let $\eta$ be its generic point and put $K:=\bC(S')=\bC(S)$. The following conditions are equivalent:
    \begin{itemize}
        \item[$(i)$] for every $s\in S$, the fibre $\cA_s$ is isomorphic, as a variety, to an abelian variety of dimension $g$.
        
        \item[$(ii)$] $\cA\to S$ is a torsor, for the \'etale  topology, under an abelian scheme over $S$ of relative dimension $g$;

        \item[$(iii)$] the generic fibre $\cA_\eta$ becomes isomorphic, as a variety, to an abelian variety of dimension $g$ after base change to an algebraic closure $\overline K$ of $K$.
    \end{itemize}
    
    The implication $(i)\Rightarrow(ii)$ is precisely the construction of the associated abelian scheme given above; see Remark~\ref{rem:localOrigins}. The implication $(ii)\Rightarrow(iii)$ follows because a torsor over $\overline{K}$ has a point and hence is trivial. For $(iii)\Rightarrow(i)$, choose an \'etale covering $\{U_i\to S\}$ by connected schemes with local sections as in Remark~\ref{rem:localOrigins}. The geometric generic fibre of each pullback is isomorphic to an abelian variety, so \cite[Theorem~6.14]{MFK} gives $\cA_{U_i}/U_i$ an abelian-scheme structure with the chosen section as identity. As the images of the $U_i$ cover $S$, condition $(i)$ follows. When these conditions hold, $f$ is a family in the sense of Definition~\ref{defn:singleFamily}, and $S$ is its smooth base.
\end{rem}

Let $f\colon\cA\to S$ be the smooth part of a family of $g$-dimensional abelian varieties, with associated abelian scheme $\cB/S$, and let $t_\sigma$ be a translation with translation vector $\sigma\colon S\to\cB$. We assume that $\bZ\sigma$ is Zariski dense in $\cB$; in particular, $\sigma$ is non-torsion and $t_\sigma$ has infinite order. Let $\Gamma_\sigma:=\langle t_\sigma\rangle
\subseteq\operatorname{Bir}(\cA')$ be the group generated by the chosen translation. We consider its action by $S$-automorphisms on the invariant open subset $\cA$. For $p\in\cA_s$, its orbit is
\[
\Gamma_\sigma(p) := \{t_\sigma^r(p):r\in\bZ\} = \{p+r\sigma(s):r\in\bZ\}.
\]
We are interested in the finite-orbit locus
\[
\fF^{(1)} := \{p\in\cA:\#\Gamma_\sigma(p)<\infty\},
\]
and ask whether it is contained in $f^{-1}(Z)$ for some
proper closed subset $Z\subsetneq S$. Since the action of $\cB_s$ on $\cA_s$ is free, for $p\in\cA_s(\bC)$ and $N\geq 1$ we have $t_\sigma^N(p)=p$ if and only if $N\sigma(s)=0_\cB(s)$. Consequently,
\begin{equation}\label{eq:singleFiniteOrbitLocus}
\fF^{(1)}=f^{-1}(T_\sigma),    
\end{equation}
where $T_\sigma$ is the set of torsion values of $\sigma$. Thus, the single-family orbit problem reduces to determining whether $T_\sigma$ is Zariski dense in $S$. This question is governed by the Betti map of the translation vector. We recall its construction and the properties needed below.

On a sufficiently small simply connected analytic open subset $\Delta\subseteq S$, choose a basis $\omega_1,\ldots,\omega_{2g}$ of the period lattice and a logarithm of $\sigma$. Write
\[
\log_{\sigma}(s)=\sum_{j=1}^{2g}\beta_j(s)\omega_j(s), \qquad \beta_j(s)\in\mathbb R.
\]
The map $\beta_{\sigma}=(\beta_1,\ldots,\beta_{2g})$ is real analytic. A change of choices replaces it by $M\beta_{\sigma}+c$, with $M\in\operatorname{GL}_{2g}(\mathbb Z)$ and $c\in\mathbb Z^{2g}$. Its generic real rank $\operatorname{rk}\beta_\sigma$ is consequently well defined. We collect some elementary properties of the Betti map in the following result.

\begin{lem}\label{lem:Betti}
The following properties hold.
\begin{itemize}
    \item[$(i)$] $\sigma(s)$ is torsion if and only if $\beta_\sigma(s)\in\mathbb Q^{2g}$, and $\operatorname{rk}\beta_\sigma\leq2\min(g,\dim S)$. Nonzero integer multiplication of $\sigma$ preserves its Betti rank.
    
    \item[$(ii)$] Betti maps commute locally with pullback, and the  Betti map of a product section is the product of the component maps.

    \item[$(iii)$] If, on a dense open subset of $S$, a nonzero multiple of $\sigma$ factors through an abelian subscheme of relative dimension $h$, then $\operatorname{rk}\beta_\sigma\leq2h$.

    \item[$(iv)$] If $\operatorname{rk}\beta_\sigma=2g$, then $T_\sigma$ is dense in $S$ for the analytic topology.
\end{itemize}
\end{lem}




Locally in the analytic topology, after choosing a basis of periods and an origin on $\cA\to S$, the translation $t_\sigma$ takes the form $x\mapsto x+\beta_\sigma(s)$ in $(\mathbb R/\mathbb Z)^{2g}$, where the Betti coordinates are taken modulo $\mathbb Z^{2g}$. Changing the origin does not change this increment: the relevant Betti map is attached to the section of $\cB/S$, not to a choice of section of $\cA\to S$.

We use the following form of relative Manin--Mumford.

\begin{thm}\label{thm:GH}
Let $\mathcal B\to S$ be an abelian scheme of relative dimension $g\geq1$ over a smooth irreducible quasi-projective complex variety, and let $\sigma \in\mathcal B(S)$ be a section such that $\mathbb Z\sigma$ is Zariski dense in $\mathcal B$. Then $T_\sigma$ is Zariski dense in $S$ if and only if $\operatorname{rk}\beta_\sigma=2g$. In particular, density of $T_\sigma$ implies $g\leq\dim S$.
\end{thm}

\begin{proof}
Apply \cite[Theorem~1.3]{GH} to the closed subvariety $\sigma(S)\subseteq\cB$. The last assertion follows from $\operatorname{rk}\beta_\sigma\leq2\dim S$.
\end{proof}

Under our assumption that $\bZ\sigma$ is Zariski dense in $\cB$, Theorem~\ref{thm:GH} gives
\begin{equation}\label{eq:singleFamilyBettiCriterion}
T_\sigma\textrm{ is Zariski dense in }S \quad\Longleftrightarrow\quad \operatorname{rk}\beta_\sigma=2g.
\end{equation}
Together with \eqref{eq:singleFiniteOrbitLocus}, \eqref{eq:singleFamilyBettiCriterion} shows that $\fF^{(1)}$ is contained in $f^{-1}(Z)$ for some proper closed subset $Z\subsetneq S$ if and only if $\operatorname{rk}\beta_\sigma<2g$. In particular, such confinement holds whenever $\dim S<g$.

\subsection{Confinement of finite orbit points on double families}\label{subsec:4.2}

We now examine a variation of the above setting. We consider two families of abelian varieties with a common total space, allowing their relative dimensions to be different and each equipped with a translation. For either translation separately, finiteness of an orbit is determined by a torsion value of its translation vector. The double-family problem concerns the additional constraints imposed by composing the two translations.

\begin{defn}\label{def:doubleFamily}
Let $g_1,g_2\geq 1$. A \emph{double family of
$(g_1,g_2)$-dimensional abelian varieties} consists of two families
\[
f_i\colon\cA'\dashrightarrow S_i',\qquad i=1,2,
\]
in the sense of Definition~\ref{defn:singleFamily}, with common nonsingular projective total space $\cA'$ and relative dimensions $g_i$. We denote their smooth parts by $f_i\colon\cA_i\longrightarrow S_i$ and their associated abelian schemes of translations by $\cB_i/S_i$. We put $d_i:=\dim S_i'=\dim S_i$.
\end{defn}

By Remark~\ref{rem:localOrigins}, each $\cA_i\to S_i$ is canonically a torsor under $\cB_i/S_i$. No origins on its fibres are chosen. Moreover, we have
\begin{equation}\label{eq:doubleDimensions}
\dim\cA'=g_1+d_1=g_2+d_2.
\end{equation}
Choose translation vectors $\sigma_i\in\cB_i(S_i)$ and denote by
\[
t_i:=t_{\sigma_i}\colon\cA_i\longrightarrow\cA_i, \qquad p\longmapsto p+\sigma_i(f_i(p)),
\]
the corresponding translations. Via \eqref{eq:translationsBirational}, we also regard each $t_i$ as a birational self-map of $\cA'$. We retain the notation
\[
T_{\sigma_i}:=\sigma_i^{-1}\bigl((\cB_i)_{\tor}\bigr) \subseteq S_i
\]
for the set of torsion values of $\sigma_i$. We shall use the following two hypotheses, which are imposed only where explicitly stated.

\begin{ass}\label{ass:1}
The rational map $\Psi:=(f_1,f_2)\colon \cA'\dashrightarrow S_1'\times S_2'$ is dominant.
\end{ass}

Hypothesis~\ref{ass:1} expresses the independence of the two fibrations: for $s$ in a dense open subset of $S_2$, the restriction
\[
\phi_s:=f_1|_{\cA_{2,s}}\colon \cA_{2,s}\dashrightarrow S_1'
\]
is dominant. Indeed, the image of $\cA_1\cap\cA_2$ under $\Psi$ contains a dense open subset $V\subseteq S_1\times S_2$, by Chevalley's theorem; see \cite[Tag~054K]{stacks-project}. For general $s\in S_2$, the fibre $V_s$ is a nonempty open subset of
$S_1$ contained in the image of $\phi_s$. The analogous assertion holds with the indices interchanged.

Moreover, Hypothesis~\ref{ass:1} implies $d_1+d_2\leq\dim\cA'$. In view of \eqref{eq:doubleDimensions}, this implies
\begin{equation}\label{eq:crossDimensionInequalities}
d_1\leq g_2,\qquad d_2\leq g_1.
\end{equation}
Under this assumption, $\Psi$ is generically finite precisely when $d_1+d_2=\dim\cA'$, or equivalently $d_1=g_2$ and $d_2=g_1$. In particular, if $g_1=g_2=g$, then $d_1=d_2=d$, and Hypothesis~\ref{ass:1} gives $d\leq g$; in this case $\Psi$ is generically finite precisely when $d=g$.

\begin{ass}\label{ass:2}
For $i=1,2$, the subset $\bZ\sigma_i$ is Zariski dense in $\cB_i$.
\end{ass}

For $i=1,2$, let $\eta_i$ be the generic point of $S_i$ and let $\overline\eta_i$ be a geometric generic point. Hypothesis~\ref{ass:2} is equivalent to $\bZ\sigma_i(\eta_i)$ being Zariski dense in $(\cB_i)_{\eta_i}$, and also to $\bZ\sigma_i(\overline\eta_i)$ being Zariski dense in $(\cB_i)_{\overline\eta_i}$. We express these equivalent conditions by saying that $\sigma_i$ generates the generic fibre. In particular, each $\sigma_i$ is non-torsion and each $t_i$ has infinite order.

The following standard observation will allow us to choose a fibre on which the translation vector generates the whole abelian variety.

\begin{lem}\label{lem:verygeneral}
Assume Hypothesis~\ref{ass:2}. For each $i=1,2$, the subgroup $\bZ\sigma_i(s)$ is Zariski dense in $(\cB_i)_s$ for every $s\in S_i(\bC)$ outside a countable union of proper closed subsets of $S_i$.
\end{lem}

\begin{proof}
Fix $i\in\{1,2\}$. Choose a countable algebraically closed subfield $F\subseteq\bC$ over which $\cB_i/S_i$ and $\sigma_i$ are defined. Write $\cB_{i,0}/S_{i,0}$ and $\sigma_{i,0}$ for models over $F$, let $\eta_{i,0}$ be the generic point of $S_{i,0}$, and put $K_0:=F(S_{i,0})$. By Hypothesis~\ref{ass:2}, the subgroup $\bZ\sigma_{i,0}(\eta_{i,0})$ is Zariski dense in $(\cB_{i,0})_{\eta_{i,0}}$.

There are only countably many proper closed $F$-subvarieties of $S_{i,0}$. Their base changes to $\bC$ are proper closed subsets of $S_i$. A point $s\in S_i(\bC)$ outside their union maps to $\eta_{i,0}$ and therefore induces an embedding $K_0\hookrightarrow\bC$. The fibre $(\cB_i)_s$ and the point $\sigma_i(s)$ are obtained from $(\cB_{i,0})_{\eta_{i,0}}$ and $\sigma_{i,0}(\eta_{i,0})$ by this extension of the ground field. Since Zariski density of rational points is preserved by field extension, $\bZ\sigma_i(s)$ is Zariski dense in $(\cB_i)_s$.
\end{proof}

We denote by
\[
\Gamma=\Gamma_{\sigma_1,\sigma_2}:=\langle t_1,t_2\rangle\subseteq\operatorname{Bir}(\cA')
\]
the group generated by the two chosen translations. A word in $t_1^{\pm1},t_2^{\pm1}$ is \emph{evaluable} at
$p\in\cA'$ if its letters, read from right to left without simplifying the expression, can be applied successively, each $t_i^{\pm1}$ being evaluated only on $\cA_i$. The empty word acts as the identity everywhere. Words representing the same birational transformation may have different evaluation domains, but their values agree wherever both are evaluable; see \cite[Tag~0A1Y]{stacks-project}.

For a collection $\Xi$ of words, set
\[
\Xi(p):=\{w(p):w\in\Xi\text{ is evaluable at }p\}, \qquad \fF(\Xi):=\{p\in\cA':\#\Xi(p)<\infty\}.
\]
In particular, $\Gamma(p)$ is the set of values of \emph{all} evaluable words in $t_1^{\pm1},t_2^{\pm1}$. Moreover, with the convention that positive powers mean repeated letters and zeroth powers mean the empty word, we consider the collection of \emph{ordered words}
\[
O_{12}:=\{t_1^n\circ t_2^r:n,r\in\bN\}
\]
and its subcollection
\[
\cT:=\{t_1^n:n\in\bN\} \cup\{t_2^r:r\in\bN\} \cup\{t_1^n\circ t_2:n\in\bN\}.
\]
Since $\Xi\subseteq\Xi'$ implies $\fF(\Xi')\subseteq\fF(\Xi)$, we have
\begin{equation}\label{eq:finiteOrbitInclusions}
\fF(\Gamma)\subseteq\fF(O_{12})\subseteq\fF(\cT).
\end{equation}

We consider the double family of Definition~\ref{def:doubleFamily}, with the notation and evaluation convention introduced above. Ordinary inverse images are taken on the domains of definition. For a closed subset $Z_i'\subseteq S_i'$, we write $f_i^{-1}[Z_i']$ for its set-theoretic total transform, as defined in our conventions; this is closed in $\cA'$ since $S_i'$ is projective. We seek to confine the finite-orbit loci to proper closed subsets of $\cA'$. For ordered orbits, a more precise question is whether $\fF(O_{12})\cap\cA_2$ lies over a proper closed subset of $S_2$, i.e.
\begin{equation}\label{eqn:orderedQuest}
    \fF(O_{12})\cap\cA_2 \subseteq f_2^{-1}(Z_2) \qquad\text{for some closed }Z_2\subsetneq S_2.
\end{equation}
For full group orbits, we ask the analogous question over both bases: whether there exist proper closed subsets $Z_i\subsetneq S_i$ such that
\[
\fF(\Gamma)\cap\cA_i\subseteq f_i^{-1}(Z_i), \qquad i=1,2.
\]

For $p\in\cA_i$ and $n\in\bZ$, every iterate $t_i^n(p)$ is defined and $t_i^n(p)=p+n\sigma_i(f_i(p))$. Since the action on the fibre is free, the set $\{t_i^n(p):n\in\bN\}$ is finite if and only if $f_i(p)\in T_{\sigma_i}$. Consider the dense open subset
\[
U_{12}:=\{p\in\cA_1\cap\cA_2:t_2(p)\in\cA_1\}.
\]
Its openness and density come from the fact that $t_2$ is an automorphism of $\cA_2$. For $p\in U_{12}$, each of the three one-parameter families that define $\cT(p)$ is fully evaluable. The preceding observation gives
\begin{equation}\label{eq:threeTorsionConditions}
    p\in\fF(\cT)
    \quad\Longleftrightarrow\quad
    \begin{cases}
        f_1(p)\in T_{\sigma_1},\\
        f_2(p)\in T_{\sigma_2},\\
        f_1(t_2(p))\in T_{\sigma_1}.
    \end{cases}
\end{equation}
Indeed, the three families are the individual translation orbits starting at $p$, and the $t_1$-orbit starting at $t_2(p)$. The additional condition $f_1(t_2(p))\in T_{\sigma_1}$ at the translated point will be called \emph{the twisted torsion condition}. It is the new constraint beyond the two single-family conditions. The locus $\fF(\cT)$ can thus be expressed as a torsion-specialization condition in an abelian scheme. Restricting all maps in these expressions to $U_{12}$, we consider the $U_{12}$-schemes
\[
\mathcal P_{12}:=f_1^*\cB_1\times_{U_{12}}f_2^*\cB_2, \qquad \mathcal P_3:=(f_1\circ t_2)^*\cB_1,
\]
and denote their natural sections by
\[
\xi_{12}:=(f_1^\ast\sigma_1,\, f_2^\ast\sigma_2), \qquad \xi_3:=(f_1\circ t_2)^\ast\sigma_1.
\]
Put
\[
\mathcal P:=\mathcal P_{12}\times_{U_{12}}\mathcal P_3,
\qquad \xi:=(\xi_{12},\xi_3).
\]
Thus~\eqref{eq:threeTorsionConditions} says that
$p\in\fF(\cT)$ if and only if $\xi(p)$ is torsion in $\mathcal P_p$, for $p\in U_{12}$.

The following lemma  shows that we can reduce our investigation to a special geometric configuration.

\begin{lem}\label{lem:reduction}
    Assume Hypotheses~\ref{ass:1} and~\ref{ass:2}. If the locus $\fF(\cT)$ is Zariski dense in $\cA'$, then
    \begin{equation}\label{eqn:reductionBettiRank}
        g_1=g_2=d_1=d_2=:g \qquad \textrm{and} \qquad \operatorname{rk}\beta_{\sigma_1} =\operatorname{rk}\beta_{\sigma_2}=2g.
    \end{equation}
\end{lem}

\begin{proof}
    Suppose that $\fF(\cT)$ is Zariski dense in $\cA'$. For each $i=1,2$, we have $\fF(\cT)\cap\cA_i\subseteq f_i^{-1}(T_{\sigma_i})$. If $\operatorname{rk}\beta_{\sigma_i}<2g_i$ for some $i$, Hypothesis~\ref{ass:2} and
    \eqref{eq:singleFamilyBettiCriterion} imply that $Z_i:=\overline{T_{\sigma_i}}^{\mathrm{Zar}} \subsetneq S_i$. Consequently,
    \[
    \fF(\cT)\subseteq (\cA'\setminus\cA_i)\cup f_i^{-1}(Z_i),
    \]
    where the set on the right is closed and proper in $\cA'$. This contradicts the assumed density. Hence $\operatorname{rk}\beta_{\sigma_i}=2g_i$ for $i=1,2$. Since $\operatorname{rk}\beta_{\sigma_i}\leq 2d_i$, we obtain $g_i\leq d_i$ for both indices. On the other hand, Hypothesis~\ref{ass:1} gives
    $d_1\leq g_2$ and $d_2\leq g_1$, by
    \eqref{eq:crossDimensionInequalities}. This proves the claim.
\end{proof}

Under Hypothesis~\ref{ass:1}, the equalities in \eqref{eqn:reductionBettiRank} imply that the first two torsion conditions in \eqref{eq:threeTorsionConditions} hold simultaneously on an analytically dense subset of $\cA_1\cap\cA_2$. More precisely,
\[
\{p\in\cA_1\cap\cA_2: f_1(p)\in T_{\sigma_1},\quad f_2(p)\in T_{\sigma_2}\}
\]
is analytically dense in $\cA_1\cap\cA_2$. Indeed, $\Psi$ is generically finite and dominant between nonsingular varieties on $\cA_1\cap\cA_2$, so in characteristic zero it is \'etale, hence locally an analytic isomorphism, on a dense open subset. Since each $T_{\sigma_i}$ is analytically dense in $S_i$, the inverse image of $T_{\sigma_1}\times T_{\sigma_2}$ is therefore analytically dense. This proves the assertion on the whole common smooth open, and hence also on $U_{12}$. Thus, the two individual torsion conditions alone cannot yield non-density in this case: the additional twisted torsion condition in \eqref{eq:threeTorsionConditions} must be analyzed. The following definition describes a relation under which the first two torsion conditions imply the third.

\begin{defn}\label{def:graphRelation}
We say that the \emph{graph relation} holds if there exist a dense open subset $U'\subseteq U_{12}$, an integer $N\geq1$, and a homomorphism of abelian schemes over $U'$,
\[
\theta\colon\mathcal P_{12}|_{U'}\longrightarrow \mathcal P_3|_{U'},
\]
such that
\begin{equation}\label{eq:mainGraphRelation}
    [N]\xi_3|_{U'}=\theta(\xi_{12}|_{U'}).
\end{equation}
\end{defn}

Let $\theta_i$ denote the restriction of $\theta$ to the $i$-th factor. The relation then takes the form
\[
[N]((f_1\circ t_2)^\ast\sigma_1) =\theta_1(f_1^\ast\sigma_1) + \theta_2(f_2^\ast\sigma_2).
\]
It is an algebraic relation between sections on a dense open subset. It implies that the first two torsion conditions force the twisted torsion condition wherever the relation is defined. The next theorem shows that the graph relation identifies the obstruction to non-density of $\fF(O_{12})$ and $\fF(\cT)$.

\begin{mainthm}\label{thm:A}
Let $f_i\colon\cA'\dashrightarrow S_i'$, $i=1,2$, be a double family of $(g_1,g_2)$-dimensional abelian varieties, with smooth parts $f_i\colon\cA_i\to S_i$ and associated abelian schemes $\cB_i/S_i$. Let $\sigma_i\in\cB_i(S_i)$ be translation vectors, and let $t_i:=t_{\sigma_i}$ be the corresponding translations. Put $d_i:=\dim S_i$. Assume that
Hypotheses~\ref{ass:1} and~\ref{ass:2} hold. Then the following conditions are equivalent.
\begin{itemize}
    \item[$(i)$] The locus $\fF(O_{12})$ is Zariski dense in $\cA'$.

    \item[$(ii)$] The locus $\fF(\cT)$ is Zariski dense in $\cA'$.

    \item[$(iii)$] There exists an integer $g\geq1$ such that
    \[
    g_1=g_2=d_1=d_2=g, \qquad \operatorname{rk}\beta_{\sigma_1} =\operatorname{rk}\beta_{\sigma_2}=2g
    \]
    and the graph relation of
    Definition~\ref{def:graphRelation} holds.
\end{itemize}
\end{mainthm}

The next theorem identifies the geometric restriction imposed by the graph relation. No equality of dimensions or assumption on the Betti ranks is required.

\begin{mainthm}
\label{thm:B}
Let $f_i\colon\cA'\dashrightarrow S_i'$, $i=1,2$, be a double family of $(g_1,g_2)$-dimensional abelian varieties, with smooth parts $f_i\colon\cA_i\to S_i$ and associated abelian schemes $\cB_i/S_i$. Let $\sigma_i\in\cB_i(S_i)$ be translation vectors, and let $t_i:=t_{\sigma_i}$ be the corresponding translations. Then the following assertions hold.
\begin{itemize}
    \item[$(i)$] If $\cB_1/S_1$ is constant, then the graph relation of Definition~\ref{def:graphRelation} holds.

    \item[$(ii)$] Assume that Hypotheses~\ref{ass:1} and~\ref{ass:2} hold. If the graph relation holds, then $\cB_1/S_1$ is isotrivial.
\end{itemize}
\end{mainthm}

The next theorem strengthens non-density to confinement over
the bases. The conclusions distinguish the ordered-orbit
locus from the full group-orbit locus.

\begin{mainthm}
\label{thm:C}
    Let $f_i\colon\cA'\dashrightarrow S_i'$, $i=1,2$, be a double family of $(g_1,g_2)$-dimensional abelian varieties, with smooth parts $f_i\colon\cA_i\to S_i$ and associated abelian schemes $\cB_i/S_i$. Let $\sigma_i\in\cB_i(S_i)$ be translation vectors, and let $t_i:=t_{\sigma_i}$ be the corresponding translations. Put $d_i:=\dim S_i$ and $\Gamma:=\langle t_1,t_2\rangle$. Assume that Hypotheses~\ref{ass:1} and~\ref{ass:2} hold.
    Suppose that at least one of condition \eqref{eqn:reductionBettiRank} and the graph relation fails. Then the following assertions hold.
    \begin{itemize}
        \item[$(i)$] The locus $\fF(\cT)$ is not Zariski dense in $\cA'$.

        \item[$(ii)$] There exists a proper closed subset $Z_2\subsetneq S_2$ such that $\fF(O_{12})\cap\cA_2\subseteq f_2^{-1}(Z_2)$.

        \item[$(iii)$] There exist proper closed subsets $W_i\subsetneq S_i$, for $i=1,2$, such that $\fF(\Gamma)\cap\cA_i\subseteq f_i^{-1}(W_i)$.
    \end{itemize}
    Moreover, there exist proper closed subsets $Z_2'\subsetneq S_2'$ and $W_i'\subsetneq S_i'$ for $i=1,2$ satisfying $Z_2\subseteq Z_2'\cap S_2$ and $W_i\subseteq W_i'\cap S_i$ such that
    \[
    \fF(O_{12})\subseteq f_2^{-1}[Z_2'], \qquad \fF(\Gamma)\subseteq f_1^{-1}[W_1']\cap f_2^{-1}[W_2'].
    \]
    In particular, all these conclusions hold whenever $\cB_1/S_1$ is non-isotrivial.
\end{mainthm}

We finally consider individual fibres on which the ordered finite-orbit locus is Zariski dense.

\begin{mainthm}
\label{thm:D}
\label{prop:exceptional-fibres}
Assume that $g_1=g_2=g\geq1$ and that $\mathbb Z\sigma_1$ is Zariski dense in $\cB_1$. Fix $b \in S_2(\mathbb C)$ such that the restriction $\phi_b:=f_1|_{\cA_{2,b}}\colon \cA_{2,b}\dashrightarrow S_1$ is dominant. Put $K_b:=\mathbb C(\cA_{2,b})$ and $G_b:=(\cB_{1,\eta_1})_{K_b}$, where $\eta_1$ is the generic point of $S_1$, and the extension $\mathbb C(S_1)\hookrightarrow K_b$ is induced by $\phi_b$. If $\fF(O_{12})\cap \cA_{2,b}$ is Zariski dense in $\cA_{2,b}$, then
\[
\sigma_2(b)\in \mathcal B_{2,b}(\mathbb C)_{\tor} \cap\operatorname{Isog}(G_b).
\]
In particular, if $\mathcal B_{2,b}$ is simple and $G_b$ is non-isotrivial, then
\[
\operatorname{ord}(\sigma_2(b)) \leq |\operatorname{Isog}(G_b)| \leq |\operatorname{Bad}(G_b)| \min_{D\in\operatorname{Bad}(G_b)} \frac{D^g}{g!}.
\]
\end{mainthm}

Theorem~\ref{thm:D} is proved in
Subsection~\ref{subsec:exceptional-fibres}. We deduce Theorem~\ref{thm:C} from
Theorems~\ref{thm:A} and~\ref{thm:B}, whose proofs are given in the next subsection.

\paragraph{Proof of Theorem~\ref{thm:C}.}
    By Theorem~\ref{thm:A}, the failure of at least one of condition~\eqref{eqn:reductionBettiRank} and the graph relation implies that $\fF(\cT)$ is not Zariski dense in $\cA'$. This proves $(i)$. Since
    \[
    \fF(\Gamma)\subseteq\fF(O_{12})\subseteq\fF(\cT),
    \]
    the other two finite-orbit loci are not Zariski dense either.

    We first observe that, for $i=1,2$, every proper closed $t_i$-invariant subset $Y\subsetneq\cA_i$ has image contained in a proper closed subset of $S_i$. Indeed, $f_i(Y)$ is closed because $f_i\colon\cA_i\to S_i$ is proper. If $f_i(Y)=S_i$, then $Y_{\overline\eta_i}$ is a nonempty closed subset of the geometric generic fibre $(\cA_i)_{\overline\eta_i}$, invariant under translation by $\sigma_i(\overline\eta_i)$. Hypothesis~\ref{ass:2} implies that the orbit of any point of this torsor under $\bZ\sigma_i(\overline\eta_i)$ is Zariski dense. Consequently, $Y_{\overline\eta_i}=(\cA_i)_{\overline\eta_i}$, which forces $Y=\cA_i$, a contradiction.

    Put $E:=\fF(O_{12})\cap\cA_2$. For $p\in E$, the entire forward $t_2$-orbit of $p$ is defined and is contained in $O_{12}(p)$, so it is finite. Moreover,
    \[
    O_{12}(t_2p)\subseteq O_{12}(p).
    \]
    Hence $t_2(E)\subseteq E$. The finiteness of the $t_2$-orbit also gives $t_2^{-1}p=t_2^{m-1}p$ for some $m\geq1$, and therefore $t_2^{-1}(E)\subseteq E$. Thus $t_2(E)=E$. Its closure $Y_{12}$ in $\cA_2$ is proper by part $(i)$ and is $t_2$-invariant because $t_2$ is an automorphism of $\cA_2$. The preceding observation shows that
    \[
    Z_2:=f_2(Y_{12})\subsetneq S_2
    \]
    is closed, proving $(ii)$. Part $(iii)$ follows from part $(ii)$ and the symmetry of the configuration.

    We now pass to the projective bases. Properness of $f_i\colon\cA_i\to S_i$ implies that its graph is closed in $\cA'\times S_i$. Since $\cA_i$ is dense in the domain of $f_i$, it follows that
    \[
    \Gamma_{f_i}\cap(\cA'\times S_i) =\{(p,f_i(p)):p\in\cA_i\}.
    \]
    The projection $\Gamma_{f_i}\to\cA'$ is surjective, being proper with dense image. Consequently, every point outside $\cA_i$ belongs to the total transform of the boundary:
    \[
    \cA'\setminus\cA_i \subseteq f_i^{-1}[S_i'\setminus S_i].
    \]
    Define
    \[
    Z_2':=\overline{Z_2}^{\,S_2'}\cup(S_2'\setminus S_2), \qquad W_i':=\overline{W_i}^{\,S_i'}\cup(S_i'\setminus S_i).
    \]
    These are proper closed subsets, with $Z_2'\cap S_2=Z_2$ and $W_i'\cap S_i=W_i$. Combining the boundary containment with $(ii)$ and $(iii)$ gives the conclusion of the theorem.

    Finally, if $\cB_1/S_1$ is non-isotrivial, then Theorem~\ref{thm:B} $(ii)$ implies that the graph relation fails. All the preceding conclusions therefore apply.
\qed

\subsection{The graph relation criterion}\label{subsec:4.3}

In this section we prove Theorem~\ref{thm:A} and ~\ref{thm:B}. We begin with a general observation on pairs of sections of abelian schemes. Let $\mathcal G/V$ be an abelian scheme over a smooth irreducible quasi-projective complex variety $V$, with relative dimension $r$. Let $\alpha:V \to \mathcal G$ be a section such that $\mathbb Z \alpha$ is Zariski dense in $\mathcal G$. Given another abelian scheme $\mathcal Q/V$ with a section $\beta:V \to \mathcal Q$, we consider the fibre product $\mathcal G \times_V \mathcal Q/V$ with the section $\zeta:=(\alpha,\beta)$. Consider the smooth closed subgroup scheme $\overline{\mathbb Z \zeta}^{\,\mathrm{Zar}}$ over $V$ and denote by $\mathcal H$ its relative identity component, it is an abelian subscheme of $\mathcal G\times_V\mathcal Q/V$. Notice that there exists an integer $m\geq1$ such that $[m]\zeta$ is a section of $\mathcal H/V$ such that $\mathbb Z ([m]\zeta)$ is Zariski dense in $\mathcal H$. Denote the projections from $\mathcal G\times_V\mathcal Q$ by $\pr_1$ and $\pr_2$, and write
\[
\rho_1:=\pr_1|_{\mathcal H}: \mathcal H\longrightarrow\mathcal G, \qquad \rho_2:=\pr_2|_{\mathcal H}: \mathcal H\longrightarrow\mathcal Q.
\]
The homomorphism $\rho_1$ is surjective, because its image is closed and contains $\mathbb Z([m]\alpha)$ which is Zariski dense in $\mathcal G$. Consequently, $\dim(\mathcal H/V)\geq r$.

\begin{lem}\label{lem:subgroup}
With the above notation, the equality $\dim(\mathcal H/V)=r$
holds if and only if there exist a dense open subset $V'\subseteq V$, an integer $N\geq1$, and a homomorphism of abelian schemes $\theta:\mathcal G|_{V'}\longrightarrow\mathcal Q|_{V'}$ such that
\begin{equation}\label{eq:subgroup-relation}
  [N](\beta|_{V'})=\theta\circ(\alpha|_{V'}).
\end{equation}
\end{lem}

\begin{proof}
Suppose first that $\dim(\mathcal H/V)=r$. Since $\rho_1:\mathcal H\to\mathcal G$ is surjective, it is an isogeny of abelian schemes. Choose an integer $e\geq1$ annihilating $\ker\rho_1$. Then $[e]_{\mathcal H}$ factors through $\rho_1$, giving a homomorphism
\[
u:\mathcal G\longrightarrow\mathcal H, \qquad u\circ\rho_1=[e]_{\mathcal H}.
\]
Since $\rho_1$ is surjective, this also implies $\rho_1\circ u=[e]_{\mathcal G}$. Recall that $[m]\zeta$ is a section of $\mathcal H$. Applying $u\circ\rho_1=[e]_{\mathcal H}$ to this section and then composing with $\rho_2$ gives
\[
[em]_{\mathcal Q}\,\beta =\rho_2\circ u\circ([m]_{\mathcal G}\,\alpha) =[m]_{\mathcal Q}\circ\rho_2\circ u\circ\alpha.
\]
Thus, setting $N:=em$ and $\theta:=[m]_{\mathcal Q}\circ\rho_2\circ u$ we obtain \eqref{eq:subgroup-relation} with $V'=V$.

Conversely, suppose that \eqref{eq:subgroup-relation} holds over a dense open subset $V'\subseteq V$. Consider the homomorphism
\[
F:(\mathcal G\times_V\mathcal Q)|_{V'} \longrightarrow\mathcal Q|_{V'}, \qquad (p,q)\longmapsto\theta(p)-[N]q.
\]
Since $F\circ\zeta|_{V'}=0$ and $F$ is a homomorphism,
its kernel contains $\mathcal H|_{V'}$. Put $\mathcal K:=\ker(\rho_1|_{V'})$. Since the first projection vanishes on $\mathcal K$, the following diagram commutes:
\[
\begin{tikzcd}
\mathcal K \arrow[r, hook] \arrow[d, "\rho_2|_{\mathcal K}"'] & \mathcal H|_{V'} \arrow[d, "F|_{\mathcal H|_{V'}}=0"] \\
\mathcal Q|_{V'} \arrow[r, "{[-N]}"'] & \mathcal Q|_{V'}.
\end{tikzcd}
\]
It follows that $[N]\circ(\rho_2|_{\mathcal K})=0$. Moreover, the second projection identifies the kernel of the first projection on the ambient product with $\mathcal Q|_{V'}$. Its restriction therefore gives a closed immersion
\[
\mathcal K\hookrightarrow\mathcal Q|_{V'}[N].
\]
Thus $\mathcal K$ is finite over $V'$. Surjectivity of $\rho_1|_{V'}$ now gives $\dim(\mathcal H/V) = \dim(\mathcal H|_{V'}/V') = r$.\\
\end{proof}

We now return to the double-family setting and retain the notation introduced in Definition~\ref{def:graphRelation} and the discussion preceding it. In the following proof of Theorem~\ref{thm:A}, Hypotheses~\ref{ass:1} and~\ref{ass:2} are assumed.

\paragraph{Proof of  Theorem~\ref{thm:A}.}
The inclusion $\cT\subseteq O_{12}$ gives $\fF(O_{12})\subseteq\fF(\cT)$, proving $(i)\Rightarrow(ii)$.

\smallskip
\noindent
\emph{$(ii)\Rightarrow(iii)$.}
Suppose that $\fF(\cT)$ is Zariski dense in $\cA'$. Then, by Lemma~\ref{lem:reduction}, we have
\[
g_1=g_2=d_1=d_2=:g, \qquad \operatorname{rk}\beta_{\sigma_1} =\operatorname{rk}\beta_{\sigma_2}=2g.
\]
We first show that $\mathbb Z\xi_{12}$ is Zariski dense in $\mathcal P_{12}$. The dominant map
\[
  \Psi=(f_1,f_2):
  \cA_1\cap\cA_2\longrightarrow S_1\times S_2
\]
is generically finite, hence \'etale on a dense open subset. With compatible local choices, the Betti map of $\xi_{12}$ is $p\longmapsto \bigl( \beta_{\sigma_1}(f_1(p)), \beta_{\sigma_2}(f_2(p)) \bigr)$. The product Betti map on $S_1\times S_2$ has rank $4g$, so pullback along the \'etale locus gives
\begin{equation}\label{eq:jointBettiRank}
  \operatorname{rk}\beta_{\xi_{12}}=4g.
\end{equation}
Denote by $\mathcal H_{12}$ the relative identity component of the Zariski closure of $\mathbb Z\xi_{12}$. If $\mathbb Z\xi_{12}$ were not Zariski dense in $\mathcal P_{12}$, then $h:=\dim(\mathcal H_{12}/U_{12})<2g$. Thus, the section $[k]\xi_{12}$ would take values in $\mathcal H_{12}$ for an integer $k\geq1$. So Lemma~\ref{lem:Betti} would give $\operatorname{rk}\beta_{\xi_{12}}\leq2h<4g$, contradicting \eqref{eq:jointBettiRank}.

Apply the construction preceding Lemma~\ref{lem:subgroup} to
\[
V=U_{12},\qquad \mathcal G=\mathcal P_{12},\qquad \mathcal Q=\mathcal P_3,\qquad \alpha=\xi_{12},\qquad \beta=\xi_3.
\]
Thus, if $\mathcal H$ denotes the relative identity component of the Zariski closure of $\mathbb Z\xi$, there is an integer $m\geq1$ such that $s:=[m]\xi$ is a section of $\mathcal H$ whose multiples are Zariski dense in $\mathcal H$. Moreover, the projection $\mathcal H\longrightarrow\mathcal P_{12}$ is surjective, and hence $\dim(\mathcal H/U_{12})\geq2g$. By \eqref{eq:threeTorsionConditions}, the section $\xi=(\xi_{12},\xi_3)$ has Zariski-dense torsion values on $U_{12}$ and the same is true of $s=[m]\xi$. Applying Theorem~\ref{thm:GH} to $s$ as a section of $\mathcal H$, we get $\dim(\mathcal H/U_{12})\leq\dim U_{12}=2g$ and consequently
\[
\dim(\mathcal H/U_{12}) =\dim(\mathcal P_{12}/U_{12}) = 2g.
\]
Lemma~\ref{lem:subgroup} now yields a dense open subset
$U'\subseteq U_{12}$, an integer $N\geq1$, and a
homomorphism of abelian schemes $\theta: \mathcal P_{12}|_{U'}\longrightarrow\mathcal P_3|_{U'}$ such that $[N](\xi_3|_{U'}) = \theta\circ(\xi_{12}|_{U'})$. This is the graph relation, proving $(iii)$.

\smallskip
\noindent
\emph{$(iii)\Rightarrow(i)$.}
Assume $(iii)$. By Lemma~\ref{lem:Betti}, each $T_{\sigma_i}$ is analytically, and hence Zariski, dense in $S_i$. Let $U'\subseteq U_{12}$, $N\geq1$ and $\theta$ give the graph relation. By Hypothesis~\ref{ass:1} and the discussion following it, there is a dense open subset $S_2^\circ\subseteq S_2$ such that, for every $b\in S_2^\circ$, the restriction
\[
\phi_b=f_1|_{A_b}:A_b\dashrightarrow S_1, \qquad A_b=f_2^{-1}(b),
\]
is dominant. Since $f_2$ is smooth, $f_2(U')$ is a dense open subset of $S_2$. Replacing $S_2^\circ$ by its intersection with $f_2(U')$, we may also assume that every such $A_b$ meets $U'$. Fix $b\in T_{\sigma_2}\cap S_2^\circ$, and put
\[
x:=\sigma_2(b), \qquad m:=\operatorname{ord}(x).
\]
Translation by $x$ defines an automorphism $T_x$ of $A_b$. Since $A_b$ is irreducible, the open subset
\[
V_b :=\bigcap_{j=0}^{m-1}T_x^{-j}(U'\cap A_b)
\]
is dense in $A_b$. Since the restriction $\phi_b|_{V_b}$ is a dominant morphism and $T_{\sigma_1}$ is Zariski dense in $S_1$, the subset
\[
D_b:=V_b\cap\phi_b^{-1}(T_{\sigma_1})
\]
is Zariski dense in $A_b$. Let $p\in D_b$ and set $p_j:=T_x^j(p)$ for $j\geq0$. Then $p_{j+m}=p_j$, and every $p_j$ belongs to $U'$. Evaluating the graph relation at $p_j$ gives
\[
[N]\sigma_1(f_1(p_{j+1})) = \theta_{p_j} \bigl(\sigma_1(f_1(p_j)),x\bigr).
\]
Since $x$ is torsion and $\sigma_1(f_1(p))$ is torsion by definition of $D_b$, this relation implies that $\sigma_1(f_1(p_j))$ is torsion for all $j \ge 0$. Each $p_j$ therefore has a finite $t_1$-orbit, all of whose iterates are defined in $\cA_1$. Every ordered word is consequently evaluable at $p$, and
\[
O_{12}(p) = \bigcup_{j=0}^{m-1} \{t_1^n(p_j):n\in\mathbb N\}
\]
is finite. Hence $D_b\subseteq\fF(O_{12})\cap A_b$. It follows that the Zariski closure of $\fF(O_{12})\cap\cA_2$ contains $A_b$ for every $b\in T_{\sigma_2}\cap S_2^\circ$. These fibres are Zariski dense in $\cA_2$, because $T_{\sigma_2}\cap S_2^\circ$ is dense in $S_2$ and $f_2$ is open and surjective. Thus $\fF(O_{12})\cap\cA_2$ is Zariski dense in $\cA'$, proving $(i)$.
\qed\\

We now derive a consequence of Corollary~\ref{thm:translation}, applied as in part $(i)$ of Remark \ref{rem:group-action-examples}, for relations between sections. We recall some standard properties of the isotrivial part introduced in our notation:
\begin{itemize}
    \item every $K$-homomorphism $u\colon G\to H$ satisfies $u(I(G))\subseteq I(H)$;

    \item for every $\mathbb C$-automorphism $\tau$ of $K$, we have $I(\tau^*G)=\tau^*I(G)$;

    \item the quotient $G/I(G)$ has no nonzero isotrivial abelian subvariety, even after a finite extension of $K$.
\end{itemize}

\begin{prop}\label{prop:quotientrigidity}
Let $A$ be a complex abelian variety, put $K=\mathbb C(A)$, and let $x\in
A(\mathbb C)$ be a point such that $\mathbb Z x$ is Zariski dense in $A$.
Let $G$ and $G'$ be abelian varieties over $K$, with $G'$ isotrivial.
Suppose that there exist points $s\in G(K)$ and $s'\in G'(K)$, an integer
$N\geq1$, and $K$-homomorphisms
\[
u\colon G\longrightarrow T_x^*G, \qquad v\colon G'\longrightarrow T_x^*G,
\]
such that
\begin{equation}\label{rig:eq:constantrelation}
[N]T_x^*s=u(s)+v(s').
\end{equation}
Then the image of $s$ in $G/I(G)$ is torsion. In particular, if $\mathbb Z
s$ is Zariski dense in $G_{\overline K}$, then $G$ is isotrivial.
\end{prop}

\begin{proof}
Put $Q:=G/I(G)$, and denote by $s_Q\in Q(K)$ the image of $s$. By the
properties of the isotrivial part, we have
\[
u(I(G))\subseteq I(T_x^*G)=T_x^*I(G).
\]
So $u$ induces a $K$-homomorphism $u_Q \colon Q\longrightarrow T_x^*Q$.
Since $G'$ is isotrivial, we also have $v(G')\subseteq I(T_x^*G)$. Thus
\eqref{rig:eq:constantrelation} gives
\begin{equation}\label{eqn:dependence}
    u_Q(s_Q)=[N]T_x^*s_Q.
\end{equation}
Let $J\subseteq Q$ be the identity component of the Zariski closure of the
cyclic subgroup $\mathbb Z s_Q$ in $Q$. The relation \eqref{eqn:dependence}
gives
\[
u_Q(J)=[N](T_x^*J)=T_x^*J.
\]
Since $J$ and $T_x^*J$ have the same dimension, the restriction
$u_Q|_J\colon J\longrightarrow T_x^*J$ is a $K$-isogeny.
Corollary~\ref{thm:translation} therefore implies that $J$ is isotrivial.
Since $Q=G/I(G)$ has no nonzero isotrivial abelian subvariety, we obtain
$J=0$. Thus the Zariski closure of $\mathbb Z s_Q$ is finite, and $s_Q$ is
torsion.

In particular, if $\mathbb Z s$ is Zariski dense in $G_{\overline K}$, then
$\mathbb Z s_Q$ is Zariski dense in $Q_{\overline K}$. Since $s_Q$ is
torsion, this forces $Q=0$. Hence $G=I(G)$ is isotrivial.
\end{proof}

We apply Proposition~\ref{prop:quotientrigidity}
to the graph relation on a suitable fibre of the
second family.

\paragraph{Proof of Theorem~\ref{thm:B}.}

\emph{(i)} Assume that $\cB_1/S_1$ is constant. Fix an identification $\cB_1=B_0\times S_1$, where $B_0$ is a complex abelian variety, and regard $\sigma_1$ as a morphism $S_1\to B_0$. Since $\cA'$ is nonsingular, the rational map
\[
  \alpha:=\sigma_1\circ f_1\colon
  \cA'\dashrightarrow B_0
\]
extends to a morphism on $\cA'$ \cite[Theorem~3.2]{MilneAV98}. We first
construct a homomorphism of abelian schemes $h\colon\cB_2\longrightarrow
B_0\times S_2$. Choose an \'etale cover $\pi\colon V\to S_2$ over which
$\cA_2$ admits an origin section $o$, as in Remark~\ref{rem:localOrigins}.
The torsor action and subtraction in $B_0$ define a morphism
$h_o\colon\cB_2\times_{S_2}V\longrightarrow B_0\times V$ whose fibre over
$v\in V(\mathbb C)$ is
\[
h_{o,v}(q)=\alpha(o(v)+q)-\alpha(o(v)), \qquad q\in\cB_{2,\pi(v)}.
\]
Since $h_o$ preserves the zero section, it is a homomorphism
\cite[Corollary~15.2(a)]{MilneAV98}. This construction is independent of the
chosen origin. Indeed, any other origin section has the form $o'=o+r$, where
$r$ is a section of $\cB_2\times_{S_2}V$, and
\[
h_{o',v}(q) = h_{o,v}(r(v)+q)-h_{o,v}(r(v)) = h_{o,v}(q).
\]
Thus the local homomorphisms agree on overlaps and descend to a map
$h\colon\cB_2\longrightarrow B_0\times S_2$ whose fibre
$h_b\colon\cB_{2,b}\to B_0$ satisfies
\[
\alpha(p+q)-\alpha(p)=h_b(q) \qquad \textrm{for every } b\in S_2(\mathbb C), \, p\in\cA_{2,b} \textrm{ and } q\in\cB_{2,b}.
\]
On $U_{12}$, the chosen trivialization identifies both $f_1^*\cB_1$ and
$\mathcal P_3$ with $B_0\times U_{12}$. Define
\[
\theta\colon\mathcal P_{12}\longrightarrow\mathcal P_3, \qquad \theta_p(a,q)=a+h_{f_2(p)}(q).
\]
Since $t_2(p)=p+\sigma_2(f_2(p))$, we have
\[
\alpha(t_2(p)) = \alpha(p)+h_{f_2(p)} \bigl(\sigma_2(f_2(p))\bigr).
\]
Under these identifications, this is precisely $\xi_3=\theta(\xi_{12})$.
Hence the graph relation holds with $N=1$.

\smallskip
\noindent

\emph{(ii)} Let $U'\subseteq U_{12}$, $N\geq1$ and $\theta$ give the graph relation, and denote by $\theta_1$ and $\theta_2$ the restrictions of $\theta$ to the two factors of $\mathcal P_{12}$. By Hypothesis~\ref{ass:1}, the restriction $f_1|_{\cA_{2,b}}$ is dominant for $b$ in a dense open subset of $S_2$. Moreover, $f_2(U')$ is dense open, since $f_2$ is smooth on $\cA_2$. Since $S_2$ is irreducible over $\mathbb C$, the complement of a countable union of proper closed subsets meets every nonempty Zariski open subset. Lemma~\ref{lem:verygeneral} therefore allows us to choose $b\in S_2(\mathbb C)$ such that
\begin{itemize}
    \item $\cA_{2,b}\cap U'\neq\varnothing$;
    \item the rational map
    $\phi:=f_1|_{\cA_{2,b}}\colon
    \cA_{2,b}\dashrightarrow S_1$
    is dominant;
    \item $\mathbb Z\sigma_2(b)$ is Zariski dense
    in $\cB_{2,b}$.
\end{itemize}
Put
\[
A:=\cA_{2,b}, \qquad B:=\cB_{2,b}, \qquad x:=\sigma_2(b).
\]
Choose a complex point of the $B$-torsor $A$ as origin, thereby identifying
$A$ with $B$. Then $t_2|_A=T_x$ is translation by $x$. Let $W\subseteq A$ be
the domain of definition of $\phi$, and consider the abelian scheme and
section
\[
\mathcal G:=\phi^*\cB_1\longrightarrow W, \qquad \widetilde s:=\phi^*\sigma_1.
\]
Put $\Omega:=A\cap U'$. By the definition of $U_{12}$, we have
$\Omega\subseteq W\cap T_x^{-1}(W)$. Thus both $\mathcal G$ and
$T_x^*\mathcal G$ are defined over $\Omega$. Since $f_2|_A$ is constant with value $b$, there are natural identifications
\[
\left.(f_1^*\cB_1)\right|_\Omega \simeq \mathcal G|_\Omega, \qquad \left.(f_2^*\cB_2)\right|_\Omega \simeq B\times\Omega, \qquad \mathcal P_3|_\Omega \simeq (T_x^*\mathcal G)|_\Omega.
\]
Under the second identification, $f_2^*\sigma_2$ becomes the constant section $x_\Omega$ with value $x$. Consequently, restricting $\theta_1$ and $\theta_2$ to $\Omega$ gives homomorphisms of abelian schemes
\[
u\colon\mathcal G|_\Omega \longrightarrow (T_x^*\mathcal G)|_\Omega, \qquad v\colon B\times\Omega \longrightarrow (T_x^*\mathcal G)|_\Omega,
\]
and the graph relation becomes
\[
[N]T_x^*\widetilde s = u(\widetilde s)+v(x_\Omega) \qquad\text{over }\Omega.
\]
We now pass to generic fibres. Put $k:=\mathbb C(S_1)$, $K:=\mathbb C(A)$ and fix the embedding $k\hookrightarrow K$ induced by $\phi$. Denoting by $\eta$ the generic point of $A$, we have
\[
G:=\mathcal G_\eta=(\cB_{1,\eta_1})_K, \qquad s:=\widetilde s(\eta)=\sigma_1(\eta_1)\in G(K).
\]
By Hypothesis~\ref{ass:2}, $s$ generates $G_{\overline K}$. The generic fibres of the preceding homomorphisms $u, v$ therefore satisfy the hypotheses of Proposition~\ref{prop:quotientrigidity}: the translation point $x$ generates $A$, and the source $B_K$ of $v$ is constant. Hence $G$ is isotrivial.

We still need to prove that $\cB_1/S_1$ is isotrivial. Choose compatible algebraic closures $\overline k\subseteq\overline K$. Since $G$ is isotrivial, there is a complex abelian variety $G_0$ and an isomorphism
\[
(\cB_{1,\eta_1})_{\overline K} \simeq (G_0)_{\overline K}.
\]
Both abelian varieties are defined over $\overline k$, and this isomorphism and its inverse are defined over a common intermediate field $\overline k\subseteq L\subseteq\overline K$ that is finitely generated over $\overline k$. Since $\overline k$ is algebraically closed, $L/\overline k$ is regular. By \cite[Corollary~15.4(b)]{MilneAV98}, both morphisms therefore descend to $\overline k$. This implies $(\cB_{1,\eta_1})_{\overline k} \simeq (G_0)_{\overline k}$, so $\cB_1/S_1$ is isotrivial.
\qed

\begin{rem}
    Only dominance of $\phi$ is used in part~$(ii)$, the extension $K/k$ need not be finite. Consequently, neither the dimension equalities nor the Betti-rank conditions in \eqref{eqn:reductionBettiRank} are required for Theorem~\ref{thm:B}.
\end{rem}

\subsection{Exceptional fibres}\label{subsec:exceptional-fibres}

We prove Theorem~\ref{thm:D}. The main step is to deduce an isogeny between $G_b$ and its translated pullback from the simultaneous torsion of two sections. The translation groups and stabilizers below are taken with respect to the canonical action of $\cB_{2,b}$ on $\cA_{2,b}$.

\begin{proof}[Proof of Theorem~\ref{thm:D}]
Put $A:=\cA_{2,b}$ and $T:=t_2|_A$, so that $T$ is translation by $\sigma_2(b)$. For every $p\in A$, all forward iterates are defined and satisfy
\[
T^m(p)=p+[m]\sigma_2(b),\qquad m\geq0.
\]
These points belong to $O_{12}(p)$. Since the torsor action is free, nonemptiness of $\fF(O_{12})\cap A$ already implies that $\sigma_2(b)$ is torsion. Choose a dense open subset $U_b\subseteq A\cap\cA_1$ on which $\phi_b$ is a morphism to $S_1$, and consider the abelian scheme and section
\[
\mathcal G_b:=\phi_b^*\cB_1\longrightarrow U_b, \qquad s:=\phi_b^*\sigma_1.
\]
The generic fibre of $\mathcal G_b$ is $G_b$. Since $\mathbb Z\sigma_1$ is Zariski dense in $\cB_1$ and $\phi_b$ is dominant, the generic value of $s$ generates $(G_b)_{\overline{K_b}}$. The generic value of $T^*s$ likewise generates $(T^*G_b)_{\overline{K_b}}$. On the dense open subset $V:=U_b\cap T^{-1}(U_b)$, consider the section
\[
\nu:=(s,T^*s) \quad\text{of}\quad \mathcal P:= \mathcal G_b|_V\times_V(T^*\mathcal G_b)|_V.
\]
For $p\in V\cap\fF(O_{12})$, both $p$ and $T(p)$ belong to $\cA_1$, and their $t_1$-orbits are contained in $O_{12}(p)$. These orbits are therefore finite, so their translation vectors $s(p)$ and $(T^*s)(p)$ are torsion. Thus $\nu$ has Zariski-dense torsion values on $V$.

As in the proof of Theorem~\ref{thm:A}, $(ii)\Rightarrow(iii)$, apply the construction preceding Lemma~\ref{lem:subgroup} to the pair $(s,T^*s)$. After shrinking $V$, this gives an abelian subscheme $\mathcal H\subseteq\mathcal P$ and an integer $N\geq1$ such that $[N]\nu$ is a section of $\mathcal H$ whose multiples are Zariski dense in $\mathcal H$. Both projections from $\mathcal H$ to the factors of $\mathcal P$ are surjective, since the generic values of $s$ and $T^*s$ generate their respective geometric generic fibres. Moreover, $[N]\nu$ has Zariski-dense torsion values on $V$, so Theorem~\ref{thm:GH} gives
\[
g\leq\dim(\mathcal H/V)\leq\dim V=g.
\]
Thus both projections induce $K_b$-isogenies on generic fibres, and hence
\[
G_b\sim_{K_b}T^*G_b.
\]
This proves $\sigma_2(b)\in\operatorname{Isog}(G_b)$.

For the final estimate, choose a complex point of $A$ as origin, identifying $A$ with $\cB_{2,b}$. Non-isotriviality of $G_b$ implies $\operatorname{Bad}(G_b)\ne\varnothing$ by Proposition~\ref{prop:Liouville-general}. If $\cB_{2,b}$ is simple, every nonzero effective divisor on $A$ is ample. Proposition~\ref{prop:translation-group-bound} therefore shows that $\operatorname{Isog}(G_b)$ is finite. Hence $\Sigma_{G_b}^0$ is trivial, and the same proposition yields the asserted estimate.
\end{proof}

The estimate depends on the bad divisors of the individual fibre. Moreover, by Proposition~\ref{inv_isog}, $\operatorname{Isog}(G_b)$ is contained in the setwise translation stabilizer of each of $\operatorname{Bad}(G_b)$ and $\operatorname{Bad}^{\infty}(G_b)$. If either stabilizer is trivial, the density assumption in Theorem~\ref{thm:D} therefore forces $\sigma_2(b)=0$.

\subsection{Isotriviality does not determine density}\label{subsec:isoDensity}
The following two double elliptic families satisfy Hypotheses~\ref{ass:1},~\ref{ass:2} and \eqref{intro:balanced}. In both cases the associated elliptic schemes are isotrivial but nonconstant, whereas the graph relation and the density of finite orbits behave differently.

\begin{expl}

\noindent\emph{(i) The graph relation holds and finite group orbits
are dense.}
Let $E$ be a complex elliptic curve, let
$q\colon E\to E/\{\pm1\}\simeq\mathbb P^1$, and let $\mathcal A'$
be the minimal resolution of
\[
(E\times E)/\langle\iota\rangle,
\qquad
\iota(a,b)=(-a,-b).
\]
The projections induce elliptic fibrations
\[
f_1([a,b])=q(a),\qquad f_2([a,b])=q(b).
\]
Their zero sections are induced by $E\times\{0\}$ and
$\{0\}\times E$, respectively, and the diagonal induces sections
$\sigma_1,\sigma_2$. Over $S:=\mathbb P^1\setminus q(E[2])$, both associated elliptic schemes become constant after the connected finite \'etale double cover $E\setminus E[2]\to S$. Their local monodromy around each omitted point is $-1$, so neither scheme is constant. The corresponding translations are induced by
\[
\widetilde t_1(a,b)=(a,b+a),\qquad
\widetilde t_2(a,b)=(a+b,b).
\]
These automorphisms commute with $\iota$, descend to the quotient,
and induce automorphisms of $\mathcal A'$.
The joint map $(f_1,f_2)$ has degree $2$. On the double covers of
the bases, both sections become the identity map of $E$; hence they
are non-torsion and both Betti maps have rank $2$.

On the double cover of $U_{12}$ induced by $E\times E$, the three
pullback elliptic schemes are constant, and the sections have fibre
coordinates
\[
\xi_{12}(a,b)=(a,b),\qquad \xi_3(a,b)=a+b.
\]
The homomorphism $(u,v)\mapsto u+v$ is equivariant for the simultaneous
sign change on the base and all three elliptic factors. It therefore
descends to
$\theta\colon\mathcal P_{12}\to\mathcal P_3$ with
$\xi_3=\theta(\xi_{12})$. Thus the graph relation holds with $N=1$. Moreover, the images of torsion pairs outside $E[2]\times E[2]$ have finite $\Gamma$-orbits: both lifted translations and their inverses preserve $E[N]\times E[N]$ for every $N\geq1$, and preserve the complement of $E[2]\times E[2]$. These images are Zariski dense. Consequently, $\mathfrak F(\Gamma)$, $\mathfrak F(O_{12})$ and $\mathfrak F(\mathcal T)$ are all Zariski dense in $\mathcal A'$.

\smallskip
\noindent\emph{(ii) The graph relation fails and finite orbits
are not dense.}
Let $\mathcal A'$ be a smooth projective model of
\[
X^\circ=\{(x,s,t)\in\mathbb A^3:s^2=x^3+t^2+1\}
\]
on which $f_1=t$ and $f_2=s$ extend to morphisms to $\mathbb P^1$. With origins at infinity, the generic fibres are
\[
\mathcal B_{1,t}:\ y^2=x^3+t^2+1, \qquad \mathcal B_{2,s}:\ v^2=u^3+s^2-1,
\]
where $(x,y)=(x,s)$ for the first fibration and $(u,v)=(-x,t)$ for the second. Fix sections $\sigma_1(t)=(-1,t)$ and $\sigma_2(s)=(1,s)$ and let $t_i$ be translation by $\sigma_i$. The map $(f_1,f_2)$ has degree $3$, so Hypothesis~\ref{ass:1}  holds.

Both elliptic schemes have $j=0$ and become constant after adjoining, respectively, $r^6=t^2+1$ and $w^6=s^2-1$. Their nontrivial local monodromy at $t=\pm i$ and $s=\pm1$ shows that they are nonconstant. After these base changes, the sections become the nonconstant maps to $E_0:Y^2=X^3+1$ given by
\[
(t,r)\longmapsto(-r^{-2},tr^{-3}), \qquad (s,w)\longmapsto(w^{-2},sw^{-3}).
\]
Thus both sections are non-torsion and have Betti rank $2$. In particular, Hypothesis~\ref{ass:2}  and \eqref{intro:balanced} hold.

We verify that the graph relation fails. Put $K=\mathbb C(\mathcal A')$ and $h=f_1\circ t_2$. The three pullback elliptic curves over $K$ have equations $Y^2=X^3+d_j$, where
\[
d_1=t^2+1,\qquad d_2=s^2-1,\qquad d_3=h^2+1.
\]
Let $C\subset X^\circ$ be the reduced irreducible curve $t=i$, and let $D$ be its strict inverse image under $t_2$. Translation is an isomorphism near the generic point of $C$, because the second fibration is smooth there. Hence $\operatorname{ord}_D(h-i)=1$. The addition law on the second elliptic fibration gives $t_2(2,0,-3i)=(0,0,i)$. At the source point, $d_1=-8$ and $d_2=-1$, so neither vanishes
identically along $D$. Consequently,
\[
\operatorname{ord}_D(d_1)=
\operatorname{ord}_D(d_2)=0,
\qquad
\operatorname{ord}_D(d_3)=1.
\]
The first two pullback curves therefore have good reduction at $D$,
whereas the third has bad reduction of type $\mathrm{II}$.
Since good reduction is invariant under isogeny, and every nonzero
homomorphism between elliptic curves is an isogeny, this implies
\[
\operatorname{Hom}_K\bigl((\mathcal P_{12})_K,
                         (\mathcal P_3)_K\bigr)=0.
\]
As $\xi_3$ is non-torsion, no identity
$[N]\xi_3=\theta(\xi_{12})$ can hold, even after shrinking $U_{12}$.
By Theorem~A, $\mathfrak F(\mathcal T)$ is not Zariski dense;
neither are $\mathfrak F(O_{12})$ and $\mathfrak F(\Gamma)$.

Thus isotriviality alone determines neither the graph relation nor
the Zariski density of the finite-orbit loci, even when both
elliptic schemes are nonconstant and both sections have maximal
Betti rank.

\end{expl}

%% file: example-paolo-francesco.tex
\def \P{{\mathbb P}}
\def \Z{{\mathbb Z}}
\def \Q{{\mathbb Q}}
\def \R{{\mathbb R}}

The purpose of this appendix to the paper by P. Dolce and F. Tropeano ``Finite translation orbits on double families of abelian varieties'' is to remark that examples of the situation studied in this paper exist in every even dimension, and to provide 
some explicit constructions, as well as indications how to prove abstract existence results in a case which has been extensively studied by geometers. 
The general framework is as 
follows. We consider an {\bf irreducible holomorphically symplectic (IHS) manifold} $X$, that is, a simply-connected manifold 
$X$ such that 
$H^0(X, \Omega^2_X)$ is one-dimensional and generated by a nowhere degenerate form $\sigma$. We can take $X$ projective, or more generally compact K\"ahler (in the situation we are looking for, projectivity shall be automatic). A typical example of such a 
manifold is a K3 surface $S$, or, more generally, the $n$-th punctual Hilbert scheme $S^{[n]}$, parameterizing subschemes of $S$ of finite length $n$. In all explicit examples, we shall be dealing with $S^{[n]}$, but the general results are valid in the general IHS context.

It is well-known that on the second cohomology $H^2(X, \Z)$ there is an integral non-degenerate quadratic form $q$, called the Beauville-Bogomolov form, which can be seen as an analogue of the intersection form on a surface. If $X\to B$ is a fibration, then the inverse image of an ample line bundle on $B$ is nef and $q$-isotropic. Conversely, a famous ``Lagrangian'', or ``hyperk\"ahler SYZ'', conjecture, checked in all known examples, in particular for $S^{[n]}$, states that if $L$ is a nef line bundle on $X$ with $q(L)=0$, 
then some power of $L$ is base-point-free, so that its sections define a fibration $\phi=\phi_L:X\to B$. Matsushita \cite{Mats} 
proved that a non-trivial fibration on an IHS manifold is equidimensional, and all smooth fibers are lagrangian tori. In particular, if $\phi$ has a section,
one obtains a family of abelian varieties on an open subset of $X$, say $\phi^0:X^0\to B^0$.

Oguiso (\cite{O}) proved that the Picard number of the generic fiber of such a fibration is always equal to one. 
In particular, the generic fiber is simple, so that the family does not have a fixed part as soon as it is not isotrivial. In fact it is easy to deduce from \cite{Bak} or \cite{AV-parab} that no finite base-change of $\phi^0$ has a fixed part unless the family is isotrivial.

By the same reason, the multiples of any non-torsion section or multisection of a family of abelian varieties arising in this way must be Zariski-dense. 

If $f$ is an automorphism of $X$ such that its action on $H^2(X,\Z)$ preserves the class of $L$ as above, then a power of $f$ preserves
the fibration $\phi_L:X\to B$ (\cite{LB}) and acts on the smooth fibers as a translation (\cite{AV-parab}). There is a way to say whether
an automorphism $\psi$ of the Neron-Severi lattice $NS(X)\subset H^2(X,\Z)$ preserving the class of $L$ comes from an actual automorphism $f:X\to X$, see {\bf ``Hodge-theoretic Torelli theorem''} by Markman, \cite{Mark}: 
it should belong to the (Hodge) monodromy group\footnote{The monodromy group is the group of 
automorphisms of $H^2(X,\Z)$ generated by all parallel transports in families, and the Hodge monodromy group is the image of its Hodge type-preserving subgroup in the group of automorphisms of the Neron-Severi lattice.}, and it should take some ample class to an ample class. The Hodge monodromy group is of finite index in the automorphism group
of $(NS(X), q)$, so replacing any $\psi$ by a power we may assume it is in there. The ample cone is governed by so-called MBM classes, a higher-dimensional 
analogue of $(-2)$-classes on K3 surfaces (\cite{AV-MBM}, \cite{AV-MorKaw}). These are primitive classes in $H^2(X, \Z)$ of bounded negative square (\cite{AV-orbits}). Inside the cone of classes of positive square in $NS(X)\otimes \R$, the ample cone is a connected component of the complement to the union of the orthogonal hyperplanes to the MBM classes of Hodge type $(1,1)$. On all known examples of IHS manifolds, in particular on $S^{[n]}$, these classes can be described explicitely. If no MBM class is orthogonal to $L$ in $(NS(X), q)$, then, up to taking a power, an automorphism of 
the lattice which fixes 
$L$ lifts to an automorphism of $X$: indeed the image of an ample class near $L$ in $NS(X)\otimes \Q$ shall be ample, so this is a consequence of Hodge-theoretic Torelli. The automorphisms preserving $L$, up to a finite index, form a free abelian group of rank $\rho-2$, where $\rho$ is the Picard number of $X$ 
(we assume here that $\rho\geq 3$, then the statement is obtained from hyperbolic geometry, see \cite{AV-parab}). If there are such MBM classes but not too many, some automorphisms may lift, see e.g. \cite{Mats-parab}: one has to further subtract from $\rho-2$ the dimension of the subspace they generate. Such automorphisms are sometimes called {\bf parabolic}.

Let us start with the following explicit example. Let $S$ be a smooth quartic surface in $\P ^3$ (it is, of course, a K3 surface). It is well-known and easy to see that $S$ can contain only finitely many (complex) lines, so if $S$ is defined over a number field, then the lines are defined over a (possibly larger) number field too. Assume $S$ contains a line $l$. Take all planes through $l$, it is a pencil of planes (they are parameterized by $\P^1$). For each such plane $P_t$, the intersection with $S$ is $l\cup C_t$, where $C_t$ is a plane cubic. This gives a fibration $\phi:S\to \P^1$ where the smooth fibers are curves of genus 1. The line $l$ induces a multisection: indeed $l$ intersects each $C_t$ in three points. So it is a trisection. 

If $S$ contains another line $l'$, which does not intersect $l$ (this is possible, e.g. on a Fermat surface, but also on others - in fact over a codimension-two subvariety of the parameter space for quartic surfaces), this gives a section of $\phi$, indeed each $P_t$ and hence each $C_t$ intersects $l'$ at one point. 
In its turn, taking the pencil of planes $P'_t$ through $l'$, we obtain another fibration of $S$, $\phi ':S\to \P^1$, with genus one fibers $C'_t$ residual to $l'$ in the intersection of $S$ and $P'_t$, a section induced by $l$, and a trisection induced by $l'$ itself. 

On the resulting abelian schemes, these trisections are non-torsion, see e.g. \cite{Hassett-budalect} where it is explained that a torsion multisection of an elliptic fibration of a K3 surface 
cannot be a rational curve. One can also choose $S$ in such a way that it contains an additional line $m$ skew 
to both $l$ and $l'$: it shall induce an additional section of both fibrations. Keeping in mind the general theory of automorphisms of IHS manifolds and 
MBM classes, one may also produce non-torsion sections on $S$ as follows. 

\medskip

{\bf Proposition 1:} If $S$ is general with the above properties, then $S$ admits an automorphism $h$ of infinite order preserving $\phi$ and acting as a 
translation along its fibers.

\medskip

{\it Proof:} For such an $S$, the lattice $NS(X)$ is of rank 3, generated by
the classes $H$ (the hyperplane section class), $l$ and $l'$, and the class $L$ of $C_t$ is $H-l$. The orthogonal complement to $L$ is generated by $L$ itself and $H-3l'$, which has square $-20$. Hence there are no MBM classes in the orthogonal complement to $L$: indeed these have square $-2$. So the result 
follows  from Hodge-theoretic Torelli.

\medskip

We derive in particular that $S$ also has a non-torsion section $h(l')$ of $\phi$.
The same applies to $\phi'$ (with $L'=H-l'$) and gives a non-torsion section $h'(l)$.

Consider now the $k$-th punctual Hilbert scheme $S^{[k]}$ of a K3 surface $S$: it parameterizes subschemes of $S$ of length $k$, e. g. $k$-ples of distinct points, or of not necessarily distinct points with some extra structure. It is often viewed as a resolution of singularities of the $k$-th symmetric power of $S$. Any fibration $g :S\to \P^1$ naturally induces the fibration 
$g^{[k]}:S^{[k]} \to \P^k=Sym^k(\P^1)$. The fiber over a point $t_1+\dots+t_k$ (where the $t_i$ are distinct points on the projective line) is just the product 
$C_{t_1}\times C_{t_2}\times \dots \times C_{t_k}$. So this is a fibration where the fibers over an open subset of the base are $k$-dimensional tori. Any section $s$ of $g$ naturally induces a section $s^{[k]}$ of $g^{[k]}$, and non-torsion induces non-torsion. 

We are now in a position to give explicit examples of the situation considered in the paper.

\medskip

{\bf Theorem:} For each $k\geq 1$ there exist algebraic varieties $X$ of dimension $2k$ with two fibrations $\phi$ and $\phi'$ from $X$ to $\P^k$,
such that $\phi$ resp. $\phi'$ induces an abelian scheme structure without a fixed part on an open subset $U$ resp. $U'$ of $X$. Each of these fibrations has an extra non-torsion section. Moreover the
multiples of these sections are Zariski-dense in $U$, $U'$.

\medskip

{\it Proof:} Take $S$ a quartic in $\P^3$ containing two skew lines $l$ and $l'$, inducing fibrations $\phi$ and $\phi'$, and
consider $\phi^{[k]}$ and $\phi'^{[k]}$ on $X=S^{[k]}$.

\medskip

Another, maybe slightly less well-known construction is as follows, see \cite{HT-abel}. Take $S$ a complete intersection of three quadrics in $\P^5$. This is again a K3 surface.
We can arrange for $S$ to contain a rational normal cubic $C$ and to contain no lines. Let $H$ be a hyperplane section divisor, then $(H-C)^2=0$, 
so curves residual to $C$ in a hyperplane section are of square zero and genus one, this gives a fibration of $S$, and $C$ induces a multisection of degree 5. Lift this fibration to $S^{[2]}$ as before, call it $\phi$.  Remark that a point of $S^{[2]}$ is either a pair of distinct points of $S$ or a point together with a tangent direction. Through each pair of points of $S$, or a point with a tangent direction, there is a unique line $l$, and it does not intersect $S$ at any extra points (indeed, since $S$ is an intersection of quadrics, the line would be contained in $S$ otherwise). The quadrics containing $S$ are parameterized by a projective plane $\P(V)$, and those among them which contain $l$, by a line in this plane, so we have a natural map from $S^{[2]}$ to the dual projective plane $\P(V^*)$, and a fiber is naturally identified to the set of lines contained in the intersection of two quadrics, known to be an abelian surface generically (when this intersection is smooth), see e.g. \cite{R}. So we have another fibration called $\phi'$. 

\medskip

{\bf Proposition 2} The curve $C^{[2]}$ viewed as a subvariety of $S^{[2]}$ induces a (possibly rational\footnote{By a rational section we mean a section defined over a dense open subset of the base.}) section of $\phi'$.

\medskip

{\it Proof:} Indeed the intersection of two sufficiently general quadrics from $\P(V)$ and the projective space $\P^3$ generated by $C$ is a union of $C$ and one of its secant lines $l$, so that $C\cap l$ gives a distinguished point in each fiber of $\phi'$.

\medskip

Note, though, that the first fibration does not have a natural section arizing from this geometric construction. However one can impose a section, e.g. by 
requiring $S$ to contain another rational normal cubic $C'$ intersecting $C$ at two points: then $C'$ induces a section of $\phi$ and $C^{[2]}$ induces a section of $\phi^{[2]}$. One may remark that there is also an abstract existence result, which follows from the Torelli theorem for K3 surfaces and Nikulin's results on lattice embedding: for any nondegenerate
even lattice $\Lambda$ of signature $(1, \rho -1)$, $\rho \leq 10$, there exists a K3 surface with Neron-Severi group $\Lambda$ (see \cite{Mor}).

Once two fibrations are constructed, the existence of parabolic automorphisms preserving each one can be deduced in the same way as in Proposition 1: indeed the description of the Neron-Severi group and of the MBM classes on $S^{[2]}$ is well-known (the latter are the classes of square $-2$ and those classes of square $-10$ which have even pairing with all other classes in $H^2(S^{[2]}, \Z)$, see \cite{HT-fourfolds} for statements, \cite{AV-K3} for an easy proof). We check the existence of a parabolic automorphism preserving $\phi$ on $S$, and of a parabolic automorphism preserving $\phi'$ on $S^{[2]}$. The details are left to the reader.

As a final remark, let us mention that many more examples can be constructed in an ``abstract'' way, by choosing a convenient lattice $\Lambda$ of low rank (but at least three), so that there is an IHS manifold of one of the four known deformation types (e.g. deformation equivalent to the Hilbert scheme of a K3 surface) $X$ with Neron-Severi lattice $\Lambda$. As the description of the MBM classes is available, by choosing the lattice carefully it is possible to arrange for two 
Beauville-Bogomolov isotropic nef classes with no, or few, orthogonal MBM classes. Since the Lagrangian conjecture is verified, this gives two lagrangian fibrations $\phi$, $\phi'$, and by Hodge-theoretic Torelli, two groups of parabolic automorphisms $P$ resp. $P'$ preserving each. One then may study the locus of points with finite orbit with respect to the group generated by some $f\in P$ and
$f'\in P'$.

Note also that IHS manifolds with two transversal lagrangian fibrations have been constructed in \cite{KV}; as the ambient space there has Picard rank two, there are no automorphisms which are interesting for us, but a suitable modification of the construction could certainly yield some. The construction of \cite{KV} is entirely based on the Torelli theorem, so is not explicit.